%% file: main.tex
\documentclass[%
  a4paper,
  onecolumn,
  algotwoe,
]{mypreprint}

\usepackage[english]{babel}

\usepackage{graphicx}

\usepackage{tikz}
  \usetikzlibrary{calc}
\usepackage{pgfplots}
\usepackage{pgfplotstable}
  \usepgfplotslibrary{fillbetween}
  \usepgfplotslibrary{units}
  \usepgfplotslibrary{external}
  \tikzset{external/system call = {%
      pdflatex \tikzexternalcheckshellescape
      -halt-on-error
      -interaction=batchmode
      -jobname "\image" "\texsource"}}
  \tikzexternaldisable
  \pgfplotsset{compat = 1.13, unit code/.code 2 args={\si{#1#2}}}
  \pgfkeys{/pgf/number format/.cd,1000 sep={\,}}

\usepackage{caption}
\usepackage{subcaption}
\usepackage{booktabs}

\usepackage{amsmath}
\usepackage{amssymb}
\usepackage{amsthm}

\usepackage{scalerel}

\usepackage{siunitx}
\usepackage{enumitem}
\usepackage{csquotes}

\newcommand{\wh}[1]{\hstretch{2.5}{\hat{\hstretch{.4}{#1}}}}

\renewcommand{\rm}[1]{\ensuremath{\mathrm{#1}}}
\newcommand{\ve}[1]{\ensuremath{\boldsymbol{#1}}}

\newcommand{\R}{\ensuremath{\mathbb{R}}}
\newcommand{\C}{\ensuremath{\mathbb{C}}}

\newcommand{\bcK}{\ensuremath{\ve{\mathcal{K}}}}
\newcommand{\bcB}{\ensuremath{\ve{\mathcal{B}}}}
\newcommand{\bcC}{\ensuremath{\ve{\mathcal{C}}}}
\newcommand{\bcS}{\ensuremath{\ve{\mathcal{S}}}}
\newcommand{\bA}{\ensuremath{\ve{A}}}
\newcommand{\bB}{\ensuremath{\ve{B}}}
\newcommand{\bC}{\ensuremath{\ve{C}}}
\newcommand{\bE}{\ensuremath{\ve{E}}}
\newcommand{\bM}{\ensuremath{\ve{M}}}
\newcommand{\bD}{\ensuremath{\ve{D}}}
\newcommand{\bK}{\ensuremath{\ve{K}}}
\newcommand{\bX}{\ensuremath{\ve{X}}}
\newcommand{\bx}{\ensuremath{\ve{x}}}
\newcommand{\bU}{\ensuremath{\ve{U}}}
\newcommand{\bu}{\ensuremath{\ve{u}}}
\newcommand{\bY}{\ensuremath{\ve{Y}}}
\newcommand{\bby}{\ensuremath{\ve{y}}}
\newcommand{\bLambda}{\ensuremath{\ve{\Lambda}}}
\newcommand{\bV}{\ensuremath{\ve{V}}}
\newcommand{\bW}{\ensuremath{\ve{W}}}
\newcommand{\bv}{\ensuremath{\ve{v}}}
\newcommand{\bw}{\ensuremath{\ve{w}}}
\newcommand{\br}{\ensuremath{\ve{r}}}
\newcommand{\bl}{\ensuremath{\ve{\ell}}}
\newcommand{\bJ}{\ensuremath{\ve{J}}}
\newcommand{\bI}{\ensuremath{\ve{I}}}
\newcommand{\bb}{\ensuremath{\ve{b}}}
\newcommand{\bc}{\ensuremath{\ve{c}}}
\newcommand{\bG}{\ensuremath{\ve{G}}}

\newcommand{\cH}{\ensuremath{\mathcal{H}}}
\newcommand{\cL}{\ensuremath{\mathcal{L}}}
\newcommand{\Htwo}{\ensuremath{\cH_{2}}}

\newcommand{\Linf}{\ensuremath{\cL_{\infty}}}

\newcommand{\bcKr}{\ensuremath{\wh{\bcK}}}
\newcommand{\bcBr}{\ensuremath{\skew1\wh{\bcB}}}
\newcommand{\bcCr}{\ensuremath{\skew2\wh{\bcC}}}
\newcommand{\Sigmar}{\ensuremath{\wh{\Sigma}}}
\newcommand{\bXr}{\ensuremath{\skew3\wh{\bX}}}
\newcommand{\bYr}{\ensuremath{\wh{\bY}}}
\newcommand{\bMr}{\ensuremath{\skew3\wh{\bM}}}
\newcommand{\bDr}{\ensuremath{\skew2\wh{\bD}}}
\newcommand{\bKr}{\ensuremath{\skew3\wh{\bK}}}
\newcommand{\bBr}{\ensuremath{\skew2\wh{\bB}}}
\newcommand{\bCr}{\ensuremath{\skew2\wh{\bC}}}

\newcommand{\IL}{\ensuremath{\mathbb{L}}}
\newcommand{\EL}{\ensuremath{\bE_{\IL}}}
\newcommand{\AL}{\ensuremath{\bA_{\IL}}}
\newcommand{\BL}{\ensuremath{\bB_{\IL}}}
\newcommand{\CL}{\ensuremath{\bC_{\IL}}}
\newcommand{\VL}{\ensuremath{\bV_{\IL}}}
\newcommand{\WL}{\ensuremath{\bW_{\IL}}}

\DeclareMathOperator{\mspan}{span}
\DeclareMathOperator{\rank}{rank}
\newcommand{\orth}{\texttt{orth}}
\newcommand{\inv}{\ensuremath{^{-1}}}
\newcommand{\trans}{\ensuremath{^{\mkern-1.5mu\mathsf{T}}}}

\newcommand{\herm}{\ensuremath{^{\mathsf{H}}}}
\newcommand{\mherm}{\ensuremath{^{-\mathsf{H}}}}
\newcommand{\abs}[1]{\ensuremath{\lvert #1 \rvert}}

\theoremstyle{plain}\newtheorem{proposition}{Proposition}

\newcommand{\bya}[1]{\ensuremath{\left( #1 \right)}}
\newcommand{\by}[1]{\ensuremath{(#1)}}
\newcommand{\bs}{\by{s}}
\newcommand{\vw}{\ensuremath{\bV, \bW}}
\newcommand{\tf}{\ensuremath{\ve{H}}}
\newcommand{\tfr}{\ensuremath{\skew4\wh{\tf}}}
\newcommand{\convinterval}{\ensuremath{\Omega}}

\newcommand{\qim}{\ensuremath{\mathfrak{i}}}
\newcommand{\qreal}[1]{\ensuremath{\operatorname{Re}\left( #1 \right)}}
\newcommand{\qimag}[1]{\ensuremath{\operatorname{Im}\left( #1 \right)}}

\DeclareMathOperator{\relerr}{relerr}
\newcommand{\sigmaplot}{\ensuremath{\lVert \tf\by{\qim\omega} \rVert_{2}}}
\newcommand{\sigmaerr}{\ensuremath{\relerr(\omega)}}

\newcommand{\irka}{\texttt{IRKA}}
\newcommand{\straika}{\texttt{StrAIKA}}
\newcommand{\tfirka}{\texttt{TF-IRKA}}
\newcommand{\sptfirka}{\texttt{SPTF-IRKA}}
\newcommand{\soirka}{\texttt{SO-IRKA}}

\definecolor{matlabblue}{HTML}{0072BD}
\definecolor{matlaborange}{HTML}{D95319}
\definecolor{matlabyellow}{HTML}{EDB120}
\definecolor{matlabpurple}{HTML}{7E2F8E}
\definecolor{matlabgreen}{HTML}{77AC30}
\definecolor{matlablightblue}{HTML}{4DBEEE}
\definecolor{matlabred}{HTML}{A2142F}

\colorlet{colA}{matlabblue}
\colorlet{colB}{matlabgreen}
\colorlet{colC}{matlabpurple}
\colorlet{colD}{matlaborange}
\colorlet{colE}{matlabred}
\colorlet{colF}{matlablightblue}

\newcommand{\mylinewidth}{1.5pt}
\newcommand{\mymarkscale}{1.5}

\pgfplotscreateplotcyclelist{qcolorlist}{%
	colA,mark=*, line width=\mylinewidth, mark options={scale=\mymarkscale,fill=white,solid}\\%
	colB,mark=triangle*,line width=\mylinewidth,mark options={scale=\mymarkscale,fill=white,solid}\\%
	colC,mark=square*,line width=\mylinewidth,mark options={scale=\mymarkscale,fill=white,solid}\\%
	colD,mark=otimes*,line width=\mylinewidth,mark options={scale=\mymarkscale,fill=white,solid}\\%
	colE,mark=diamond*,line width=\mylinewidth,mark options={scale=\mymarkscale,fill=white,solid}\\%
	TUMGrau,mark=pentagon*,line width=\mylinewidth,mark options={scale=\mymarkscale,fill=white,solid}\\%
	colA,mark=otimes*,line width=\mylinewidth,mark options={scale=\mymarkscale,fill=white,solid}\\%
	colB,mark=diamond*,line width=\mylinewidth,mark options={scale=\mymarkscale,fill=white,solid}\\%
	colC,mark=pentagon*,line width=\mylinewidth,mark options={scale=\mymarkscale,fill=white,solid}\\%
	colD,mark=*,line width=\mylinewidth,mark options={scale=\mymarkscale,fill=white,solid}\\%
	colE,mark=triangle*,line width=\mylinewidth,mark options={scale=\mymarkscale,fill=white,solid}\\%
	colF,mark=square*,line width=\mylinewidth,mark options={scale=\mymarkscale,fill=white,solid}\\%
}

\pgfplotsset{colormap={qmap}{color=(colA); rgb255=(229,229,229); color=(colC)}}

\newcommand{\plotinterval}[2]{%
	\addplot [name path=min, mark=none, black!40, forget plot] coordinates {(#1, \pgfkeysvalueof{/pgfplots/ymin}) (#1, \pgfkeysvalueof{/pgfplots/ymax})};%
	\addplot [name path=max, mark=none, black!40, forget plot] coordinates {(#2, \pgfkeysvalueof{/pgfplots/ymin}) (#2, \pgfkeysvalueof{/pgfplots/ymax})};%
	\addplot[fill=black!5, area legend] fill between[of=min and max];%
}

\pgfplotstableset{
  columns/algorithm/.style={column name = Algorithm, string type, column type={l}},
  columns/linf_error/.style={column name = $\relerr_{\Linf, \Omega}$, string type, column type={S[round-mode=places, round-precision=2,scientific-notation = true, table-format=1.2e-1]}},
  columns/n_iter/.style={column name = $n_\rm{iter}$, string type, column type={S[fixed-exponent=0, table-omit-exponent, table-format=3]}},
  columns/n_ls/.style={column name = $n_\rm{ls}$, string type, column type={S[fixed-exponent=0, table-omit-exponent, table-format=4]}},
  columns/t_c/.style={column name = $t_\rm{c} \left[\si{\second}\right]$, string type, column type={S[fixed-exponent=0, table-omit-exponent, round-mode=places, round-precision=2,table-format=4.2]}},
  columns/mark_maxiter/.style={column name = ,string type, column type={l}},
}

\usepackage{cleveref}

\begin{document}

\title{Adaptive choice of near-optimal expansion points for interpolation-based
  structure-preserving model reduction}

\author[$\ast$]{Quirin Aumann}
\affil[$\ast$]{Max Planck Institute for Dynamics of Complex Technical Systems,
  Sandtorstra{\ss}e 1, 39106 Magdeburg, Germany.\authorcr
  \email{aumann@mpi-magdeburg.mpg.de }, \orcid{0000-0001-7942-5703}}

\author[$\dagger$]{Steffen W. R. Werner}
\affil[$\dagger$]{Courant Institute of Mathematical Sciences,
  New York University, New York, NY 10012, USA.
  \email{steffen.werner@nyu.edu}, \orcid{0000-0003-1667-4862}}

\shorttitle{Near-optimal structure-preserving model reduction}
\shortauthor{Q. Aumann, S.~W.~R. Werner}
\shortdate{}
\shortinstitute{}



\abstract{%
    Interpolation-based methods are well-established and effective approaches 
    for the efficient generation of accurate reduced-order surrogate models. 
    Common challenges for such methods are the automatic selection of good or 
    even optimal interpolation points and the appropriate size of the 
    reduced-order model.
    An approach that addresses the first problem for linear, unstructured
    systems is the Iterative Rational Krylov Algorithm  (IRKA), which computes
    optimal interpolation points through iterative  updates by solving linear
    eigenvalue problems.
    However, in the case of preserving internal system structures, optimal
    interpolation points are unknown, and heuristics based on nonlinear
    eigenvalue problems result in numbers of potential interpolation points
    that typically exceed the reasonable size of reduced-order systems.
    In our work, we propose a projection-based iterative interpolation method
    inspired by IRKA for generally structured systems to adaptively compute
    near-optimal interpolation points as well as an appropriate size for the
    reduced-order system.
    Additionally, the iterative updates of the interpolation points can be
    chosen such that the reduced-order model provides an accurate approximation
    in specified frequency ranges of interest.
    For such applications, our new approach outperforms the established methods
    in terms of accuracy and computational effort.
    We show this in numerical examples with different structures.
}

\novelty{%
  We propose a new algorithm for the construction of interpolating, structured,
  reduced-order models via projection.
  The method efficiently determines new interpolation points from
  the solutions of low-dimensional linear eigenvalue problems, adaptively
  chooses an appropriate size of the reduced-order model, and can be used to
  obtain high-fidelity approximations in limited frequency ranges of interest.
}

\maketitle


\section{Introduction}%
\label{sec:intro}

Simulation, control and optimization of dynamical systems are essential for
many applications.
In this work, we consider structured linear systems in the frequency (Laplace)
domain of the form
\begin{equation} \label{eqn:freqsys}
  \Sigma\colon \begin{cases}
    \bcK\bs \bX\bs = \bcB\bs \bU\bs, \quad
    \bY\bs = \bcC\bs \bX\bs,
  \end{cases}
\end{equation}
where $\bcK\colon \C \to \C^{n \times n}$ describes the system dynamics,
$\bcB\colon \C \to \C^{n \times m}$ the system's input and
$\bcC\colon \C \to \C^{p \times n}$ the system's output behavior;
see~\cite{BeaG09} for motivation of~\cref{eqn:freqsys} example systems.
The functions $\bX\colon \C \to \C^{n}$,
$\bU\colon \C \to \C^{m}$ and
$\bY\colon \C \to \C^{p}$ denote the internal states,
inputs and outputs, respectively.
For all $s \in \C$ for which $\bcK$ is invertible, and $\bcB$ and $\bcC$ can be
evaluated, the corresponding transfer function
$\tf\colon \C \to \C^{p \times m}$ directly
relates the system's inputs to outputs:
\begin{equation} \label{eqn:transfun}
  \tf\bs = \bcC\bs \bcK\bs\inv \bcB\bs.
\end{equation}

The most commonly considered  structure of dynamical systems is given by
first-order differential equations in the time domain
\begin{equation} \label{eqn:fosys}
  \Sigma\colon \begin{cases}
    \bE \dot{\bx}\by{t} = \bA \bx\by{t} + \bB \bu\by{t}, \quad
    \bby\by{t} = \bC \bx\by{t},
  \end{cases}
\end{equation}
with the system matrices $\bA, \bE \in \R^{n \times n}$,
$\bB \in \R^{n \times m}$ and $\bC \in \R^{p \times n}$.
Systems of the form~\cref{eqn:fosys} are also commonly referred to as
\emph{unstructured systems} due to them being considered as the standard case.
Applying the Laplace transformation to~\cref{eqn:fosys} yields a
frequency-dependent system of the form~\cref{eqn:freqsys}, with
\begin{equation*}
   \bcK\bs = (s \bE - \bA)\inv, \quad
   \bcB\bs = \bB, \quad
   \bcC\bs = \bC,
\end{equation*}
and the corresponding transfer function $\tf\bs = \bC (s \bE - \bA)\inv \bB$.
On the other hand, the modeling of specific physical phenomena hands down other
differential structures into dynamical systems.
The modeling of mechanical structures, structural vibrations, wave movement or
electrical circuits classically leads to differential equations with
second-order time derivatives, which in frequency domain are described by
transfer functions of the form
\begin{equation} \label{eqn:sotf}
  \tf\bs = (\bC_{\rm{p}} + s\bC_{\rm{v}}) (s^{2} \bM + s \bD + \bK)\inv \bB;
\end{equation}
see, for example,~\cite{AbrM87, Wer21, Wu16} and references therein.
A different structure occurs when modeling incomplete systems resulting in
time-delay structures, which are expressed as exponential terms in the
frequency domain, e.g., with the transfer function
\begin{equation} \label{eqn:tdtf}
  \tf\bs = \bC(s \bE - \bA_{\rm{0}} - \operatorname{e}^{-\tau s}
    \bA_{\rm{d}})\inv \bB,
\end{equation}
for some constant time delay $\tau > 0$; see, e.g.,~\cite{GaoK19}.
Many other structures exist in literature, which are used to model,
for example,
poroelasticity~\cite{AumDJetal22},
viscoelasticity~\cite{VanMM13},
or interior acoustic problems~\cite{CohHKetal08}.
While some structures, like second-order systems~\cref{eqn:sotf}, can be
reformulated into standard form~\cref{eqn:fosys}, this is not necessarily
possible for all occurring structures, including time-delay
systems~\cref{eqn:tdtf}.
Also, by reformulation into unstructured form, the number of states increases
and structure inherent properties are typically lost in subsequent
computational procedures.

In general, there is a demand for highly accurate models in practical
applications.
As a result, the number of equations describing~\cref{eqn:freqsys} vastly grows
and the computational efficient solution of~\cref{eqn:freqsys} in terms of
resources such as time and memory is often impossible.
Model order reduction methods are a remedy to this problem as they aim for the
construction of cheap-to-evaluate yet accurate surrogate models that approximate
the systems' input-to-output behavior while being described by a significantly
smaller number of equations $r \ll n$, which eases the demand on computational
resources required for the evaluation of the systems.
Many model reduction techniques have been developed for unstructured
systems~\cref{eqn:fosys}; see, for example,~\cite{Ant05}.
In addition, the preservation of internal system structures such
as~\cref{eqn:sotf,eqn:tdtf} is desired as this typically yields more accurate
approximations as well as the preservation of structure inherent properties.
Also, if the reduced-order model is to be coupled to other systems, preserving
the structure is advantageous because the same coupling conditions as for the 
full-order model can be applied to the reduced surrogate~\cite{DecDMetal21}.

Several structure-preserving model order reduction methods have been developed
in recent years.
Many of these have been tailored to particular structures that occur, for
example, in vibrational problems~\cite{HetTF12, AumW23, BedBDetal20, Wer21},
network systems~\cite{CheKS17, EggKLetal18},
or systems with Hamiltonian structure~\cite{BenZ22, GugPBetal12, HesPR22}.
The framework in~\cite{BeaG09} allows the reduction of dynamical systems
with arbitrary internal structures based on transfer function interpolation.
The quality of reduced-order models obtained by interpolation strongly depends
on the choice of interpolation points.
Therefore, a variety of strategies has been developed to perform successive
greedy searches for suitable interpolation points based on estimating the
approximation
error~\cite{CheFB22, BonFAetal16, FenKB15, RumGD14, PanWL13, AumM23, FenLBetal22}
or computing the exact error in, for example, the
$\Linf$-norm~\cite{AliBMetal20, SchV21}.

On the other hand, the \emph{Iterative Rational Krylov Algorithm}
(\irka{}) is a well established interpolation method for unstructured
systems~\cref{eqn:fosys} that iteratively updates the interpolation
points~\cite{GugAB08}.
At convergence, the interpolating reduced-order model satisfies the necessary
$\Htwo$-optimality conditions.
Several extensions of \irka{} for structured systems using similar ideas have
been proposed.
For second-order systems~\cref{eqn:sotf}, the \soirka{} method
from~\cite{Wya12} aims for an iterative process similar to \irka{}.
In~\cite{AumM22}, this has been considered as basis for a method to choose the
resulting approximation order adaptively.
\emph{Transfer Function IRKA} (\tfirka{})~\cite{BeaG12} can be applied
to arbitrarily structured systems and yields $\Htwo$-optimal but unstructured
reduced-order models.
A structure-preserving variant of \tfirka{} has been proposed in~\cite{SinGB16}.

A different take on structured model order reduction are data-driven methods.
Since here only measurements of the transfer function~\cref{eqn:transfun} are
used to compute realizations of dynamical systems, the original structure
can be arbitrary.
One of the most well-known approaches of this type is the Loewner framework,
which constructs a reduced-order model that interpolates
provided data samples~\cite{MayA07}.
The original formulation of the Loewner framework only considers the
construction of unstructured systems~\cref{eqn:fosys}, but it has been
extended to find structured realizations in~\cite{SchUBetal18}.
Recently, structured extensions of the barycentric form for second-order
systems~\cref{eqn:sotf} have been proposed that allow the extension of further
data-driven frequency domain methods to the structure-preserving
setting~\cite{GosGW23, WerGG22}.

In this work, we present a new approach to compute accurate reduced-order models
that preserve the internal structure of the original system.
Based on an \irka{}-like iteration scheme, the new method computes in
every step a new set of interpolation points (and tangential directions)
which are then employed in the structure-preserving interpolation
framework~\cite{BeaG09}.
Instead of considering nonlinear eigenvalue problems corresponding to the
resolvent terms of the structured systems, the Loewner framework allows us
to solve instead linear eigenvalue problems in each step and to determine the
approximation order adaptively and with respect to limited frequency
regions of interest if desired.

The remainder of this manuscript is structured as follows:
After introducing the mathematical preliminaries in \Cref{sec:math}, we revisit
the structure-preserving transfer function \irka{} and extend that method to
the case of multiple-input/multiple-output systems in \Cref{sec:sptfirka}.
Our new model reduction method is then described in \Cref{sec:straika}.
In \Cref{sec:numerics}, a number of numerical experiments are used to compare
the new method to established model reduction techniques.
The paper is concluded in \Cref{sec:conclusions}.


\section{Mathematical preliminaries}
\label{sec:math}


\subsection{Structure-preserving interpolation via projection}%
\label{sec:structpres}

We consider here interpolation-based model order reduction methods,
which compute surrogate models approximating the dynamics of the high-fidelity
system~\cref{eqn:freqsys} while having much smaller dimensions $r \ll n$.
Structure-preserving model order reduction methods construct approximations
of~\cref{eqn:freqsys} with the same internal structure
\begin{equation} \label{eqn:freqsys_red}
  \Sigmar\colon \begin{cases}
  \bcKr\bs \bXr\bs = \bcBr\bs \bU\bs, \quad
  \bYr\bs = \bcCr\bs \bXr\bs,
  \end{cases}
\end{equation}
where $\bcKr\colon \C \to \C^{r \times r}$,
$\bcBr\colon \C \to \C^{r \times m}$,
$\bcCr\colon \C \to \C^{p \times r}$ and
$\bXr\colon \C \to \C^{r}$,
$\bYr\colon \C \to \C^{p}$.
Additionally, the compositions of the matrix-valued functions
in~\cref{eqn:freqsys,eqn:freqsys_red} are the same:
If the center term in~\cref{eqn:freqsys} is given in frequency-affine form
\begin{equation} \label{eqn:freqaff}
  \bcK\bs = \sum\limits_{j = 1}^{n_{\bcK}} g_{j}\bs \bcK_{j},
\end{equation}
with $g_{j}\colon \C \to \C$ and constant matrices
$\bcK_{j} \in \C^{n \times n}$, for $j = 1, \ldots, n_{\bcK}$, then
the center term of the reduced-order model must have the form
\begin{equation} \label{eqn:freqaff_red}
  \bcKr\bs = \sum\limits_{j = 1}^{n_{\bcK}} g_{j}\bs \bcKr_{j},
\end{equation}
where $\bcKr_{j} \in \C^{r \times r}$, for $j = 1, \ldots, n_{\bcK}$.
Since the scalar functions in~\cref{eqn:freqaff,eqn:freqaff_red} are identical,
the internal system structure is preserved and the system matrices $\bcK_{j}$
of the full-order system can be replaced by their reduced-order counterparts
$\bcKr_{j}$.
The same relations must hold for the input and output terms $\bcB$
and $\bcC$.
Consider as an example the second-order system with transfer
function~\cref{eqn:sotf} from \Cref{sec:intro}.
A structure-preserving reduced-order model will be of the form
\begin{equation*}
  \tfr\bs = (\bCr_{\rm{p}} + s\bCr_{\rm{v}}) (s^{2} \bMr + s \bDr + \bKr)\inv
    \bBr.
\end{equation*}

To act as suitable surrogate, the reduced-order model must approximate the
input-to-output behavior of the original system at least for a range of
frequencies $s \in \C$, which are important for the application in question.
In other words, the outputs of~\cref{eqn:freqsys,eqn:freqsys_red}
should match up to a specified tolerance $\tau$ in appropriate norms for a
given input:
\begin{equation*}
    \lVert \bY - \bYr \rVert \leq \tau \cdot \lVert \bU \rVert.
\end{equation*}
For many time domain models and, in general, frequency domain models, the
relation above can be reformulated in terms of the transfer functions of the
original and reduced-order model such that
\begin{equation*}
    \lVert \tf - \tfr \rVert \leq \tau
\end{equation*}
holds.

Following~\cite{BeaG09}, any matrix-valued function of the
form~\cref{eqn:transfun} can be interpolated by a reduced-order transfer
function $\tfr$, while preserving the internal system structure using the
projection approach.
Given two reduction spaces with basis matrices
$\vw \in \C^{n \times r}$, the reduced-order model is computed by
\begin{equation} \label{eqn:projection}
  \bcKr\bs = \bW\herm \bcK\bs \bV, \quad 
  \bcBr\bs = \bW\herm \bcB\bs, \quad
  \bcCr\bs = \bcC\bs \bV.
\end{equation}
While there are many potential choices for the basis matrices $\bV$ and $\bW$,
we concentrate here on transfer function interpolation, i.e., the
matrices $\vw$ are constructed such that the transfer function $\tfr$
corresponding to~\cref{eqn:projection} interpolates the full-order transfer
function~\cref{eqn:transfun} at chosen points.
The following proposition gives a concise overview.

\begin{proposition}[Structured {interpolation~\cite[Thm.~1]{BeaG09}}]%
  \label{prp:strint}
  Let $\tf$ be the transfer function~\cref{eqn:transfun} of a linear system,
  described by~\cref{eqn:freqsys}, and $\tfr$ the reduced-order transfer
  function constructed via projection~\cref{eqn:projection}.
  Let the matrix functions $\bcC$, $\bcK\inv$, $\bcB$ and $\bcKr\inv$
  be analytic in the interpolation point $\sigma \in \C$.
  Then, the following statements hold.
  \begin{enumerate}[label = (\alph*)]
    \item If $\mspan \left( \bcK\by{\sigma}\inv \bcB\by{\sigma} \right)
      \subseteq \mspan(\bV)$ holds, then $\tf\by{\sigma} = \tfr\by{\sigma}$.
    \item If $\mspan \left( \bcK\by{\sigma}\mherm \bcC\by{\sigma}\herm \right)
      \subseteq \mspan(\bW)$ holds, then $\tf\by{\sigma} = \tfr\by{\sigma}$.
    \item If $\bV$ and $\bW$ are constructed as above,
      then additionally $\tf{}'\by{\sigma} = \tfr{}'\by{\sigma}$ holds.
	\end{enumerate}
\end{proposition}

Overall, only linear systems of equations need to be solved for the
construction of the basis matrices $\bV$ and $\bW$ in \Cref{prp:strint}.
However, the interpolation point $\sigma$ has to be known beforehand and its
choice has a large influence on the approximation quality of the resulting
reduced-order model.
Traditionally, points are chosen linearly or logarithmically equidistant on the
frequency axis $\qim\R$ in frequency ranges of interest to reduce the worst
case approximation error of the transfer function given by the $\Linf$-norm.
This typically leads to a reasonable approximation behavior over the considered
frequency range but easily misses features of the system, which are not close
enough to the interpolation points, or may result in unnecessarily large
reduced-order models.


\subsection{Unstructured interpolation via the Loewner framework}%
\label{sec:loewner}

Independent of the structure of the original system, the Loewner framework
can be used to construct unstructured systems from transfer function
evaluations~\cite{MayA07, AntLI17}.
As we will use this framework at several points throughout this manuscript,
it is summarized below following the description in~\cite{AntLI17}.

Given $2 q$ transfer function measurements
$\tf_{k} := \tf\by{s_{k}} \in \C^{p \times m}$ at some locations $s_{k} \in \C$,
for $k = 1, \dots, 2 q$, the data is partitioned into two sets
\begin{equation*}
  \begin{cases}
    \bya{\kappa_{i}, \br_{i}, \bw_{i}},~\text{where}~
      \kappa_{i} = s_{i},~
      \bw_{i} = \tf_{i} \br_{i},~\text{for}~i = 1, \dots, q,\\
    \bya{\mu_{j}, \bl_{j}, \bv_{j}},~\text{where}~
      \mu_{j} = s_{q + j},~
      \bv_{j}\herm = \bl_{j}\herm \tf_{q + j},
      ~\text{for}~j = 1, \dots, q,
  \end{cases}
\end{equation*}
with right and left tangential directions $\br_{i} \in \C^{m}$ and
$\bl_{j} \in \C^{p}$, for $i, j = 1, \ldots, q$.
In practice, it has been shown to be beneficial for numerical reasons to
partition the data in an alternating way with respect to the ordering of the
absolute values of the sampling points.
Under the assumption that the sets of sampling points are disjunct,
$\{ \kappa_{i} \}_{i = 1}^{q} \cap \{ \mu_{j} \}_{j = 1}^{q} = \emptyset$,
the partitioned data is arranged in the Loewner and shifted Loewner
matrices
\begin{align} \label{eqn:loewL}
  \IL & =
    \begin{bmatrix}
      \displaystyle \frac{\bv_{1} \br_{1} - \bl_{1} \bw_{1}}
      {\mu_{1} - \kappa_{1}}
      & \cdots &
      \displaystyle \frac{\bv_{1} \br_{q} - \bl_{1} \bw_{q}}
      {\mu_{1} - \kappa_{q}} \\
      \vdots & \ddots & \vdots \\
	   \displaystyle \frac{\bv_{q} \br_{1} - \bl_{q} \bw_{1}}
      {\mu_{q} - \kappa_{1}}
      & \cdots &
      \displaystyle \frac{\bv_{q} \br_{q} - \bl_{q} \bw_{q}}
      {\mu_{q} - \kappa_{q}}
    \end{bmatrix}, \\ \label{eqn:loewLs}
  \IL_{\sigma} & =
    \begin{bmatrix}
      \displaystyle \frac{\mu_{1} \bv_{1} \br_{1} -
      \kappa_{1} \bl_{1} \bw_{1}}
      {\mu_{1} - \kappa_{1}}
      & \cdots &
      \displaystyle \frac{\mu_{1} \bv_{1} \br_{q} -
      \kappa_{q} \bl_{1} \bw_{q}}
      {\mu_{1} - \kappa_{q}} \\
      \vdots & \ddots & \vdots \\
      \displaystyle \frac{\mu_{q} \bv_{q} \br_{1} -
      \kappa_{1} \bl_{q} \bw_{1}}
      {\mu_{q} - \kappa_{1}}
      & \cdots &
      \displaystyle \frac{\mu_{q} \bv_{q} \br_{q} -
      \kappa_{q} \bl_{q} \bw_{q}}
      {\mu_{q} - \kappa_{q}}
	\end{bmatrix}.
\end{align}
If the matrix pencil $\IL_{\sigma} - \lambda \IL$ is regular, i.e., there exist
a $\lambda \in \C$ such that
$\det\left(\IL_{\sigma} - \lambda \IL \right) \neq 0$, and given 
the matrices
\begin{equation} \label{eqn:loewWV}
  \WL = \begin{bmatrix} \bw_{1} & \ldots & \bw_{q} \end{bmatrix} \quad
  \text{and} \quad
  \VL = \begin{bmatrix} \bv_{1}\herm \\ \vdots \\ \bv_{q}\herm \end{bmatrix},
\end{equation}
the transfer function of the form
$\tf_{\IL}\bs = \WL (\IL_{\sigma} - s \IL)\inv \VL$
tangentially interpolates the given data such that
$\tf_{\IL}\by{\kappa_{i}} \br_{i} = \bw_{i}$ and
$\bl_{j}\herm \tf_{\IL}\by{\mu_{j}} = \bv_{j}\herm$ hold,
for $i, j = 1, \ldots, q$.
The corresponding state-space realization of the underlying dynamical system
is then given by
\begin{equation} \label{eqn:loew}
  \EL := -\IL, \quad
  \AL := -\IL_{\sigma}, \quad
  \BL := \VL, \quad
  \CL := \WL,
\end{equation}
using the matrices from~\cref{eqn:loewL,eqn:loewLs,eqn:loewWV}.

The rank of the Loewner pencil
$n_{\IL} = \rank \bya{\begin{bmatrix} \IL & \IL_{\sigma} \end{bmatrix}}$
directly uncovers the minimal order of a model, which is required to
interpolate the given data.
In practice, it is reasonable to truncate any redundant data, which might have
been collected into $(-\IL, -\IL_{\sigma}, \WL, \VL)$.
The required truncation matrices $\bV$ and $\bW$ can be chosen as the right and
left singular vectors obtained from singular value decompositions (SVDs) of
$\begin{bmatrix} \IL & \IL_{\sigma} \end{bmatrix}$ and
$\begin{bmatrix} \IL\herm & \IL_{\sigma}\herm \end{bmatrix}\herm$.
Truncating the matrices of singular vectors at $n_{\IL}$ columns and
projecting the Loewner realization~\cref{eqn:loew} yields a model interpolating
the given data.
Truncating after $r_{\IL} < n_{\IL}$ results in a model approximating the
provided data.

In general, without further modifications, the models obtained from the Loewner
framework may have complex-valued matrices.
However, many systems are described by real-valued matrices in practical
applications.
Under the assumption that the original transfer function follows the
\emph{reflection principle}, i.e., $\overline{\tf\bs} = \tf\by{\overline{s}}$
holds for all $s \in \C$ for which $\tf$ is defined, sampling points as well
as transfer function data and tangential direction can be chosen closed
under conjugation, i.e., if $\kappa$ is a sampling point so is
$\overline{\kappa}$, and $\tf\by{\kappa}$ and
$\tf\by{\overline{\kappa}} = \overline{\tf\by{\kappa}}$ are the corresponding
complex conjugate transfer function measurements.
In this case, there exists a state-space transformation for~\cref{eqn:loew}
that yields real-valued matrices.
Assuming that all given data is complex-valued, closed under conjugation,
and ordered into complex pairs, then the transformation to obtain real-valued
matrices is given by
\begin{equation*}
  \bJ = \bI_{q} \otimes \left(\frac{1}{\sqrt{2}} \begin{bmatrix}
    1 & -\qim \\ 1 & \qim \end{bmatrix} \right),
\end{equation*}
and the transformed system
$(-\bJ \herm \IL \bJ,\ -\bJ\herm \IL_{\sigma} \bJ,\ \bJ\herm  \WL,\ \VL \bJ)$
has real-valued matrices and satisfies the same interpolation conditions as the
original Loewner system~\cref{eqn:loew}; see~\cite{AntLI17}.


\subsection{Unstructured \texorpdfstring{$\Htwo$}{H2}-optimal interpolation}
\label{sec:tfirka}

\begin{algorithm}[t]
  \SetAlgoHangIndent{1pt}
  \DontPrintSemicolon
  \caption{Transfer function IRKA (\tfirka{}).}%
  \label{alg:tfirka}
  
  \KwIn{Transfer function $\tf\bs$,
    initial interpolation points $\left\{ \sigma_{j} \right\}_{j = 1}^{r}$ and
    tangential directions $\left\{ \bb_{j} \right\}_{j = 1}^{r}$ and
    $\left\{ \bc_{j} \right\}_{j = 1}^{r}$.}
  \KwOut{Reduced first-order model
    $\Sigmar\colon \left( \EL, \AL, \BL, \CL \right)$.}
  
  \While{no convergence}{
    Construct $\EL$, $\AL$, $\BL$ and $\CL$ as
      in~\cref{eqn:tfirka_A,eqn:tfirka_BC,eqn:tfirka_E}
      using the interpolation points $\left\{ \sigma_{j} \right\}_{j = 1}^{r}$,
      and tangential directions $\left\{ \bb_{j} \right\}_{j = 1}^{r}$ and
      $\left\{ \bc_{j} \right\}_{j = 1}^{r}$.\;
    
    Compute the generalized eigenvalues and eigenvectors
      $\left\{\left( \lambda_{j}, \bx_{j}, \bby_{j} \right)
      \right\}_{j = 1}^{r}$ from
      \begin{equation*}
        \AL \bx_{j} = \lambda_{j} \EL \bx_{j} \quad \text{and} \quad
        \bby_{j}\herm \AL = \lambda_{j} \bby_{j}\herm \EL.
      \end{equation*}\;\vspace{-\baselineskip}

    Update the interpolation points and tangential directions via
      \begin{equation*}
        \sigma_{j} \gets -\lambda_{j}, \quad
        \bb_{j}\herm \gets \bby_{j}\herm \BL \quad \text{and} \quad
        \bc_{j} \gets \CL \bx_{j},
      \end{equation*}
      for $j = 1, \ldots, r$.\;
  }
\end{algorithm}

A different approach for the construction of reduced-order models for structured
systems is the \emph{Transfer Function IRKA} (\tfirka{}) from~\cite{BeaG12}.
Like the original \irka{}~\cite{GugAB08}, this method computes $\Htwo$-optimal
approximations but can also be applied to structured systems
like~\cref{eqn:freqsys}, because only evaluations of the transfer function and
its derivative are needed for the algorithm.
However, \tfirka{} computes a reduced-order model of first-order
form~\cref{eqn:fosys}, i.e., the approach can be applied to structured systems
but does not preserve the structure in the reduced-order model.

The procedure of \tfirka{} is as follows:
Instead of computing an interpolating realization of the reduced-order model by
projection, an interpolating first-order realization is obtained using the
Loewner framework~\cite{MayA07} in every iteration step.
In contrast to the variant of the Loewner framework described in the
previous section, the two sets of interpolation points are chosen to be
identical.
This leads to a modification of the formulas~\cref{eqn:loewL,eqn:loewLs} 
involving the derivative of the sampled transfer function.
Given a transfer function $\tf\bs$, its derivative $\tf{}'\bs$,
interpolation points $\left\{ \sigma_{j} \right\}_{j = 1}^{r}$,
and right and left tangential directions
$\left\{ \bb_{j} \right\}_{j = 1}^{r}$ and
$\left\{ \bc_{j} \right\}_{j = 1}^{r}$, with $\bb_{j} \in \C^{m}$ and
$\bc_{j} \in \C^{p}$, this variant of the Loewner framework
constructs a first-order model
$\tf_{\IL}\bs = \CL \left(s \EL - \AL \right)\inv \BL$ satisfying the following
tangential Hermite interpolation conditions:
\begin{equation*}
  \tf\by{\sigma_{j}} \bb_{j} = \tf_{\IL}\by{\sigma_{j}} \bb_{j}, \quad
  \bc_{j}\herm \tf\by{\sigma_{j}} =
    \bc_{j}\herm \tf_{\IL}\by{\sigma_i}, \quad
  \bc_{j}\herm \tf{}'\by{\sigma_{j}} \bb_{j} =
    \bc_{j}\herm \tf_{\IL}'\by{\sigma_{j}} \bb_{j},
\end{equation*}
for all $j = 1, \ldots, r$.
The entries of the matrices in the Loewner realization are constructed via
\begingroup
\allowdisplaybreaks
\begin{align}
  \label{eqn:tfirka_E}
  \left( \EL \right)_{i, j} & :=
    \begin{cases}
      \displaystyle -\frac{\bc_{i}\herm \big(\tf\by{\sigma_{i}} -
        \tf\by{\sigma_{j}} \big) \bb_{j}}{\sigma_{i} - \sigma_{j}} &
        \text{if}~ i \neq j, \\
      \displaystyle -\bc_{i}\herm \tf{}'\by{\sigma_{i}} \bb_{i} &
        \text{if}~ i = j,
    \end{cases} \\
  \label{eqn:tfirka_A}
  \left( \AL \right)_{i, j} & :=
    \begin{cases}
      \displaystyle -\frac{\bc_{i}\herm \big( \sigma_{i} \tf\by{\sigma_{i}}
        - \sigma_{j} \tf\by{\sigma_{j}} \big) \bb_{j}}
        {\sigma_{i} - \sigma_{j}} &
        \text{if}~ i \neq j, \\
      \displaystyle -\bc_{i}\herm \big( s \tf\by{s} \big){}'(\sigma_{i})
        \bb_{i} &
        \text{if}~ i = j,
    \end{cases}\\
  \label{eqn:tfirka_BC}
  \BL & := \begin{bmatrix} \bc_{1}\herm \tf\by{\sigma_1} \\ \vdots \\
    \bc_{r}\herm \tf\by{\sigma_{r}} \end{bmatrix}
  \quad \text{and} \quad
  \CL := \begin{bmatrix} \tf\by{\sigma_{i}} \bb_{1} & \ldots &
    \tf\by{\sigma_{r}} \bb_{r} \end{bmatrix}.
\end{align}
\endgroup

Similar to the classical \irka{} method, the eigenvectors and mirror images of
the eigenvalues of $\AL -\lambda \EL$ with respect to the imaginary axis are
used as interpolation points and tangential directions in the next iteration
step of \tfirka{}.
At convergence, the algorithm yields a reduced-order model with a first-order
state-space realization satisfying the first-order interpolatory
$\Htwo$-optimality conditions~\cite{GugAB08,BeaG12}.
Note that in the case that the high-dimensional system also has first-order
structure, \irka{} and \tfirka{} are equivalent and converge to the same
reduced-order model~\cite{BeaG12}.
The main steps of \tfirka{} are summarized in \Cref{alg:tfirka}.
Since the reduced-order model is directly obtained from the underlying Loewner
framework, the realness of the original model can be preserved using the
technique described in \Cref{sec:loewner}.


\section{Structure-preserving near-optimal interpolation}%
\label{sec:spnoi}

In the following, we consider two iteration schemes similar to \irka{} for
finding near-optimal interpolation points for structure-preserving model
reduction in the case of general systems~\cref{eqn:freqsys} with transfer
functions of the form~\cref{eqn:transfun}.
Before we derive our new approach in \Cref{sec:straika}, we generalize the 
$\Htwo$-norm based method from~\cite{SinGB16} to the case of
multiple-input/multiple-output (MIMO) systems.
As it follows similar concepts, we use this method as the main benchmark for
the performance of our new approach in the numerical experiments.


\subsection{Structure-preserving transfer function IRKA}
\label{sec:sptfirka}

\begin{algorithm}[t]
  \SetAlgoHangIndent{1pt}
  \DontPrintSemicolon
  \caption{Structure-preserving transfer function IRKA (\sptfirka{}).}%
  \label{alg:sptfirka}
  
  \KwIn{Dynamical system $\Sigma\colon (\bcK, \bcB, \bcC)$,
    initial interpolation points $\left\{ \sigma_{j} \right\}_{j = 1}^{r}$ and
    tangential directions $\left\{ \bb_{j} \right\}_{j = 1}^{r}$ and
    $\left\{ \bc_{j} \right\}_{j = 1}^{r}$.}
  \KwOut{Reduced-order system
    $\Sigmar\colon (\bcKr, \bcBr, \bcCr)$.}
  
  \While{no convergence}{
    Compute the orthogonal basis matrices \vw{} via
      \begin{align*}
		\bV & \gets \orth \left( \begin{bmatrix}
          \bcK(\sigma_{1})\inv \bcB\by{\sigma_{1}} \bb_{1} & \ldots &
          \bcK\by{\sigma_{r}}\inv \bcB\by{\sigma_{r}} \bb_{r} \end{bmatrix}
          \right),\\
		\bW & \gets \orth \left( \begin{bmatrix}
          \bcK\by{\sigma_{1}}\mherm \bcC\by{\sigma_{1}}\herm \bc_{1}
          & \ldots &
          \bcK\by{\sigma_{r}}\mherm \bcC\by{\sigma_{r}}\herm \bc_{r}
          \end{bmatrix} \right).
      \end{align*}\;\vspace{-\baselineskip}
      \label{alg:sptfirka_bases}
  
    Project the system matrices such that
      \begin{equation*}
        \bcKr\bs \gets \bW\herm \bcK\bs \bV, \quad
        \bcBr\bs \gets \bW\herm \bcB\bs, \quad
        \bcCr\bs \gets \bcC\bs \bV.
      \end{equation*}\;\vspace{-\baselineskip}
    
    Compute an order-$r$ approximation
      \begin{equation*}
        \bcS_{r}\bs = \CL \left(s \EL - \AL \right)\inv \BL,
      \end{equation*}
      by applying \tfirka{} (\Cref{alg:tfirka}) to
      $\tfr\bs = \bcCr\bs \bcKr\bs\inv \bcBr\bs$.\;
      \label{alg:sptfirka_tfirka}

    Compute the generalized eigenvalues and eigenvectors
      $\left\{\left( \lambda_{j}, \bx_{j}, \bby_{j} \right)
      \right\}_{j = 1}^{r}$ from
      \begin{equation*}
        \AL \bx_{j} = \lambda_{j} \EL \bx_{j} \quad \text{and} \quad
        \bby_{j}\herm \AL = \lambda_{j} \bby_{j}\herm \EL.
      \end{equation*}\;\vspace{-\baselineskip}
      \label{alg:sptfirka_ev}
   
    Update the interpolation points and tangential directions via
      \begin{equation*}
        \sigma_{j} \gets -\lambda_{j}, \quad
        \bb_{j}\herm \gets \bby_{j}\herm \BL \quad \text{and} \quad
        \bc_{j} \gets \CL \bx_{j},
      \end{equation*}
      for $j = 1, \ldots, r$.\;
      \label{alg:sptfirka_upd}
  }
\end{algorithm}

The problem of constructing near-optimal interpolants for general structured
systems has been considered before in~\cite{SinGB16}.
Therein, the authors present an $\Htwo$-norm inspired strategy based on
\tfirka{} in combination with the structured interpolation framework from
\Cref{prp:strint} to compute structure-preserving reduced-order models.
\sptfirka{}, as sketched in \Cref{alg:sptfirka}, can in general be seen as a
two-step approach:
First, a structured reduced-order model is computed via projection using
\Cref{prp:strint};
then, the transfer function of this structured reduced-order model is
approximated by \tfirka{}, which yields an $\Htwo$-optimal first-order
realization, from which the mirror images of its poles are then used to update
the interpolation points for the next iteration.
The resulting reduced-order model is structure-preserving due to the employed
projection framework; the unstructured realization obtained from
\tfirka{} is only used to update the interpolation points.

Originally, \sptfirka{} has been formulated for
single-input/single-output (SISO) systems in~\cite{SinGB16}.
The extension to the MIMO case in \Cref{alg:sptfirka} follows directly from
the observation that \tfirka{} yields tangential interpolation
conditions for MIMO systems.
Consequently, basis matrices ensuring tangential interpolation are constructed
in \Cref{alg:sptfirka_bases} of \Cref{alg:sptfirka}.
Similar to the interpolation points, the tangential directions are updated
in every step of the iteration in \Cref{alg:sptfirka_upd} of
\Cref{alg:sptfirka} by computing additionally to the eigenvalues also the
corresponding left and right eigenvectors of the $\Htwo$-optimal
approximation in \Cref{alg:sptfirka_ev}.
Note that also the tangential version of \tfirka{} is used in
\Cref{alg:sptfirka_tfirka} as given in \Cref{alg:tfirka}.
Realness of reduced-order models computed with \sptfirka{} can be
preserved similarly to the procedure used in the original \irka{} and stated
in~\cite[Cor.~2.2]{GugAB08}:
Given a set of interpolation points with tangential directions, which is closed
under complex conjugation, the basis matrices \vw{} can be chosen to be
real valued.
Using these matrices in a projection~\cref{eqn:projection} preserves the
realness of the original system matrices, while enforcing the tangential
interpolation conditions.

In terms of computational effort for the choice of $\Htwo$-optimal
interpolation points, the two-step approach in \Cref{alg:sptfirka} can be seen
as beneficial.
Applying \tfirka{} (\Cref{alg:tfirka}) directly to the full-order system
requires the solution of $r$ linear systems of equations of order $n$ in each
step of the iteration.
In \sptfirka{} however, \tfirka{} is only applied to transfer functions
for which linear systems of dimension $r$ have to be solved.
In this situation, the outer loop of \sptfirka{} (\Cref{alg:sptfirka}) can also
be seen as a pre-reduction step that reduces the computational costs of
\tfirka{}.
Similar ideas to reduce the computational costs of iterative model-order
reduction methods have been used, for example,
in~\cite{AumM22, CasL18, BenGH15}.

Besides losing $\Htwo$-optimality in \sptfirka{}, an important difference
between \tfirka{} and \sptfirka{} lies in the requirements of the methods on
the availability of the original system.
\tfirka{} is a true black-box approach, where only access to transfer function
evaluations are required.
In contrary, \sptfirka{} requires access to the system matrices to construct
the basis matrices \vw{} as well as for the projection step.


\subsection{Structure-preserving adaptive iterative Krylov algorithm}
\label{sec:straika}

\begin{algorithm}
  \SetAlgoHangIndent{1pt}
  \DontPrintSemicolon
  \caption{Structure-preserving adaptive iterative Krylov algorithm
    (\straika{}).\hspace*{-\baselineskip}}%
  \label{alg:straika}

  \KwIn{Dynamical system $\Sigma\colon (\bcK, \bcB, \bcC)$,
    initial interpolation points $\left\{ \sigma_{j} \right\}_{j = 1}^{r}$,
    tangential directions $\left\{ \bb_{j} \right\}_{j = 1}^{r}$ and
    $\left\{ \bc_{j} \right\}_{j = 1}^{r}$,
    frequency range $\convinterval$,
    Loewner sampling points $\left\{ \theta_{j} \right\}_{j = 1}^{2q}$,
    maximum reduced order $r_{\max}$.}
  \KwOut{Reduced-order system
    $\Sigmar\colon (\bcKr, \bcBr, \bcCr )$ of order $r \leq r_{\max}$.}

  \While{no convergence}{
    Compute the orthogonal basis matrices \vw{} via
      \begin{align*}
		\bV & \gets \orth \left( \begin{bmatrix}
          \bcK\by{\sigma_{1}}\inv \bcB\by{\sigma_{1}} \bb_{1} & \ldots &
          \bcK\by{\sigma_{r}}\inv \bcB\by{\sigma_{r}} \bb_{r} \end{bmatrix}
          \right),\\
		\bW & \gets \orth \left( \begin{bmatrix}
          \bcK\by{\sigma_{1}}\mherm \bcC\by{\sigma_{1}}\herm \bc_{1}
          & \ldots &
          \bcK\by{\sigma_{r}}\mherm \bcC\by{\sigma_{r}}\herm \bc_{r}
          \end{bmatrix} \right).
      \end{align*}\;\vspace{-\baselineskip}
      \label{alg:straika_bases}
  
    Project the system matrices such that
      \begin{equation*}
        \bcKr\bs \gets \bW\herm \bcK\bs \bV, \quad
        \bcBr\bs \gets \bW\herm \bcB\bs, \quad
        \bcCr\bs \gets \bcC\bs \bV.
      \end{equation*}\;\vspace{-\baselineskip}
      \label{alg:straika_proj}
	
    Compute the Loewner interpolant~\cref{eqn:loew},
      $\Sigma_{\IL} \colon \left(\EL, \AL, \BL, \CL\right)$,
      via~\cref{eqn:loewL,eqn:loewLs,eqn:loewWV} with samples of $\tfr\bs$ in
      the points $\left\{ \theta_{j} \right\}_{j = 1}^{2q}$ and random
      left and right tangential directions.\;
      \label{alg:straika_loew}
	
    Compute the generalized eigenvalues and eigenvectors
      $\left\{\left( \lambda_{j}, \bx_{j}, \bby_{j} \right)
      \right\}_{j = 1}^{k}$ from
      \begin{equation*}
        \AL \bx_{j} = \lambda_{j} \EL \bx_{j} \quad \text{and} \quad
        \bby_{j}\herm \AL = \lambda_{j} \bby_{j}\herm \EL.
      \end{equation*}\;\vspace{-\baselineskip}
      \label{alg:straika_eig}
		
	Choose eigentriples according to frequency region of interest 
      \begin{equation*}
        \bLambda_{\Theta} \gets \left\{
          \left(\lambda_{j}, \bx_{j}, \bby_{j} \right)
          \Big\vert~
          \abs{\qimag{\lambda_{j}}} \subset \convinterval,
          ~\text{for}~ j = 1, \ldots, k \right\}.
      \end{equation*}\;\vspace{-\baselineskip}
      \label{alg:straika_pselect}
 
	\eIf{$\abs{\bLambda_{\Theta}} > r_{\max}$}{ \label{alg:dominance}
      Compute dominance~\cref{eqn:dom} for all poles in
        $\bLambda_{\Theta}$ with respect to $\Sigma_{\IL}$.\;

      Keep only the $r = r_{\max}$ most dominant poles in
        $\bLambda_{\Theta}$.\;
	}{
      Set $r \gets \abs{\bLambda_{\Theta}}$.
    }

    Update the interpolation points and tangential directions via
      \begin{equation*}
        \sigma_{j} \gets -\lambda_{j}, \quad
        \bb_{j}\herm \gets \bby_{j}\herm \BL \quad \text{and} \quad
        \bc_{j} \gets \CL \bx_{j},
      \end{equation*}
      such that
      $(\lambda_{j}, \bx_{j}, \bby_{j}) \in \bLambda_{\Theta}$,
      for $j = 1, \ldots, r$.\;
  }
\end{algorithm}

In addition to accuracy problems already observed in the original
publication~\cite{SinGB16}, a flexible application of \sptfirka{} and \tfirka{}
is limited by the fact that the final reduced order $r$ has to be fixed before
the algorithm is started.
The choice of a reasonable $r$ is highly problem dependent and an a priori
choice can often only be based on heuristics or in-depth knowledge about the
system dynamics.
In the cases where a maximum $r$ is not given by implementational restrictions,
it needs to be determined by several independent runs of \tfirka{} or
\sptfirka{} followed by system evaluations to estimate the approximation errors.
Another limitation of many \irka{}-like methods is that the user has no
influence on the distribution of the interpolation points.
For some applications, surrogates that approximate the high-fidelity model in a
specific frequency range only are more interesting than global approximations.
While frequency-limited variants of \irka{} exist~\cite{VuiPA13}, these methods
rely on the first-order realization of the full- as well as the reduced-order
model and are computationally costly for large-scale systems.

Here, we present a new approach for the structure-preserving realization of
reduced-order models building on similar concepts as \sptfirka{}, but also
addressing the issues raised above.
We may call this new method the
\emph{\textbf{Str}ucture-preserving \textbf{A}daptive \textbf{I}terative
\textbf{K}rylov \textbf{A}lgorithm} (\straika{}).
The approach is summarized in \Cref{alg:straika}.


\subsubsection{Computational procedure}

Comparing \Cref{alg:sptfirka,alg:straika}, the main computational procedures
look similar.
Struc\-ture-preserving reduced-order models are computed via \Cref{prp:strint},
which are then used to construct Loewner surrogates that are used to
update the interpolation points and tangential directions for the next
iteration step.
The main difference between \sptfirka{} and \straika{} lies in the construction
of the Loewner interpolants during the iteration.
While \sptfirka{} employs a complete run of \tfirka{} to construct an order-$r$
unstructured approximation of the structured reduced-order model
$\Sigmar$, in \sptfirka{} the transfer function $\tfr$ is sampled in the
frequency range of interest $\convinterval$ to reveal all essential system
dynamics.
Similar to the methods discussed in~\cite{CasL18, AumM22} the intermediate
models $\Sigmar$ and $\Sigma_{\IL}$ are used to leverage the computational
costs of different tasks.

In \Cref{alg:straika_loew} of \Cref{alg:straika}, we use the variant of the 
Loewner framework described in \Cref{sec:loewner} that only relies on the
evaluation of the low-dimensional transfer function $\tfr$ rather than its
derivatives as needed in \tfirka{}.
However, this can be arbitrarily replaced by other approaches for the
identification of unstructured first-order systems~\cref{eqn:fosys} from
frequency domain data.
This includes other variants of the Loewner framework such as the one described
in \Cref{sec:tfirka}, its block version~\cite{MayA07}, and variations in
these for choosing the dominant dynamics~\cite{KarGA21}, but also completely
different methods can be employed such as vector fitting~\cite{GusS99, DrmGB15},
RKFIT~\cite{BerG17} or the AAA algorithm~\cite{NakST18}.
The additional computational cost of evaluating $\tfr$, which is of order
$r \leq r_{\max}$, is negligible compared to updating the basis matrices $\vw$,
which requires decompositions of large-scale matrices of dimension $n$, if
cases with $r_{\max} \ll n$ are considered.
The advantages of considering \Cref{alg:straika_loew} detached from the
desired reduced order are that, first, concepts such as oversampling and
localized sampling can be used to influence the accuracy of the approximation
$\Sigma_{\IL}$ in the frequency range $\Omega$ of interest, and second,
that the amount of poles in the frequency range of interest $\Omega$ is a
strong indicator for the reduced order needed to well approximate the original
transfer function in this region.

Realness of the original system matrices can be preserved in the reduced-order
model throughout the iteration using similar ideas as for the previously
discussed methods.
Under the assumption that the initial interpolation points and tangential
directions are closed under complex conjugation, and the original model has a
reflective transfer function, real-valued matrices $\vw$ can be computed in
\Cref{alg:straika_bases} of \Cref{alg:straika} by splitting basis contributions
corresponding to complex conjugate interpolation points and by concatenating
\begin{equation*}
  \bV_{\R} = \begin{bmatrix} \qreal{\bV} & \qimag{\bV} \end{bmatrix} \quad
  \text{and} \quad
  \bW_{\R} = \begin{bmatrix} \qreal{\bW} & \qimag{\bW} \end{bmatrix}.
\end{equation*}
Thereby, the matrices of the reduced-order model computed in
\Cref{alg:straika_proj} are also real-valued.
Using the realification of the Loewner framework as described at the end
of \Cref{sec:loewner} leads to sets of eigenvalues and eigenvectors in
\Cref{alg:straika_eig} closed under conjugation.


\subsubsection{Interpolation point selection}

\begin{figure}[t]
  \centering
  \begin{subfigure}[t]{\textwidth}
    \centering
  \tikzexternalenable%
  \tikzsetnextfilename{region_spectrum}%
  \input{graphics/region_spectrum.tikz}%
  \tikzexternaldisable%

    \caption{Spectrum of first-order realization with selected interpolation
      points.}
    \label{fig:region_spectrum}
  \end{subfigure}

  \vspace{.5\baselineskip}
  \begin{subfigure}[b]{.49\textwidth}
    \centering
  \tikzexternalenable%
  \tikzsetnextfilename{region_tf}%
  \input{graphics/region_tf.tikz}%
  \tikzexternaldisable%

    \caption{Transfer functions.}
    \label{fig:region_tf}
  \end{subfigure}%
  \hfill%
  \begin{subfigure}[b]{.49\textwidth}
    \centering
  \tikzexternalenable%
  \tikzsetnextfilename{region_err}%
  \input{graphics/region_err.tikz}%
  \tikzexternaldisable%

    \caption{Relative approximation error.}
    \label{fig:region_err}
  \end{subfigure}

  \vspace{.5\baselineskip}
  \tikzexternalenable%
  \tikzsetnextfilename{region_legend}%
  \input{graphics/region_legend.tikz}%
  \tikzexternaldisable%

  \caption{Approximation of a model in a specified frequency region
    \convinterval.
    Only mirror images with respect to the imaginary axis of poles in the
    specified frequency region are considered leading to an accurate local
    approximation of the transfer function.}
  \label{fig:region}
\end{figure}

Reduced-order models computed with the structure-preserving framework presented
in \Cref{sec:structpres} approximate the full-order model well in the vicinity
of chosen interpolation points.
This observation can be used to compute reduced-order models, which approximate
the original model in a specific frequency region only.
Additionally, \irka{}-like methods aim for interpolation at the mirror images
of system poles with respect to the imaginary axis.
\Cref{fig:region} illustrates the combination of these two ideas.
The part of the spectrum of the model close to the imaginary axis is shown in
\Cref{fig:region_spectrum} and all eigenvalues $\lambda_{j}$ with
$3000 < \abs{\qimag{\lambda_{j}}} < 4000$ are marked, which corresponds to the
frequency region of interest $\Omega = [3000, 4000]$\,\SI{}{\radian\per\s}.
The mirror images of these eigenvalues are considered as interpolation points
in \straika{} for the structure-preserving interpolation framework.
The transfer function of the interpolating structured reduced-order
model is shown in \Cref{fig:region_tf} with the pointwise relative approximation
error in \Cref{fig:region_err}.
It can be seen that the reduced-order model is an accurate approximation of the
original system in the vicinity of the interpolation points $\sigma_{j}$.
Choosing only interpolation points in the frequency region $\convinterval$,
which is important for the application of the reduced-order model, can
therefore be a strategy to decrease the required size of the reduced-order
model.

Note that there are no limitations for the choice of $\Omega$.
For the global approach, it can be chosen to be the complete positive real axis
$\R_{\geq 0}$ but in other applications, only subintervals might be of interest
such that 
\begin{equation*}
  \Omega = \bigcup\limits_{j = 1}^{k} [\omega_{1, j}, \omega_{2, j}],
\end{equation*}
where $\omega_{1, j}, \omega_{2, j} \in \R_{\geq 0}$ and
$\omega_{1, j} \leq \omega_{2, j}$.
In \Cref{alg:straika_pselect} of \Cref{alg:straika}, the absolute value of the
imaginary part of the eigenvalues is considered, which implies a certain
symmetry in importance of positive and negative frequency regions; see also
\Cref{fig:region_spectrum}.
For certain applications, it may be advantageous to select eigenvalues using
other criteria, for example, their distance to the imaginary axis.

In principle, the reduced-order model might grow too large if all interpolation
points inside a defined region are considered.
Especially in cases where the global dynamics are approximated, the order $r$
can grow fast and even approach $n$.
For such cases, \straika{} chooses eigenvalues up to a defined maximum
$r_\rm{max}$ as locations for interpolation points. 
In this case, the $r_\rm{max}$ interpolation points, which will result in a
suitably good approximation of the original model, have to be selected from all
potential interpolation points.
To this end, the \emph{dominance} of all poles given by the eigenvalues in
\convinterval{} is computed and, for the $r_\rm{max}$ most dominant poles, the
interpolation points for the next iteration are set as their mirror images.
The systems constructed in \Cref{alg:straika_loew} of \Cref{alg:straika} are in
first-order unstructured form~\cref{eqn:fosys}.
For such systems, the dominance of a pole $\lambda_{j}$ with corresponding
right and left eigenvectors $\bx_{j}$ and $\bby_{j}$ is defined as
\begin{equation} \label{eqn:dom}
    d_{j} = \frac{\left\lVert \lambda_{j} (\bC_{\IL} \bby_{j})
      (\bx_{j} \herm \bB_{\IL}) \right\rVert_{2}}{\abs{\qreal{\lambda_{j}}}};
\end{equation}
where $\bC_{\IL}$ and $\bB_{\IL}$ are the output and input matrices
constructed in \Cref{alg:straika_loew} of \Cref{alg:straika}, respectively.
A pole $\lambda_{j}$ is called dominant, if $d_{j} > d_{k}$ for all $j \neq k$;
see~\cite{MarLP96}.

To ensure a high approximation quality in the frequency range of interest
$\convinterval$, it is often beneficial to include also the first potential
interpolation points located outside both ends of $\convinterval$.
While having only a small impact on the size of the reduced-order model, this 
can greatly increase the accuracy of the reduced-order model, especially, if
the full-order model has poles near the boundaries of $\convinterval$.


\section{Numerical experiments}
\label{sec:numerics}

We now demonstrate the performance of \straika{} in comparison with the
established \irka{}-like methods \tfirka{} and \sptfirka{}.
Where applicable, the classical \irka{} is also included in the comparison. 
The numerical experiments have been performed on a laptop equipped with an
AMD Ryzen\texttrademark~7~PRO~5850U and \SI{12}{\giga\byte} RAM running on 
Linux Mint~21 as operating system.
All algorithms have been implemented and run with
MATLAB\textsuperscript{\textregistered} version 9.11.0.1837725
(R2021b Update~2).
The results for \irka{} have been computed with M-M.E.S.S.
version~2.2~\cite{SaaKB22}.
The source code, data and results of the numerical experiments are available
at~\cite{supAumW23}.

In all experiments, we use a maximum number of iterations:
$n_{\rm{iter, max}} = 50$ for \straika{} and the outer iterations of
\sptfirka{}, and $n_{\rm{iter, max}} = 100$ for \tfirka{} and the inner
iterations of \sptfirka{}.
The algorithms terminate, if the relative difference between the interpolation
points in two consecutive iterations falls under the threshold of
$\epsilon = \num{1e-3}$.
To compare the accuracy of the methods, we plot the pointwise relative
approximation errors, given by
\begin{equation*}
  \sigmaerr := \frac{\lVert \tf\by{\qim \omega} - \tfr\by{\qim \omega}
    \Vert_{2}}{\lVert \tf\by{\qim \omega}\rVert_{2}}
\end{equation*}
in specified frequency intervals of interest
$\omega \in [\omega_{\min}, \omega_{\max}]$.
We also  we approximate the local, relative errors in $\convinterval$ under
the $\Linf$-norm via
\begin{equation*}
  \relerr_{\Linf, \Omega} = \frac{\max\limits_{\omega \in \convinterval}
    \lVert \tf\by{\qim \omega} - \tfr\by{\qim \omega} \lVert_{2}}
    {\max\limits_{\omega \in \convinterval}
    \lVert \tf\by{\qim \omega} \lVert_{2}}
    \approx \frac{\lVert \tf - \tfr \rVert_{\Linf, \Omega}}
    {\lVert \tf \rVert_{\Linf, \Omega}},
\end{equation*}
using equidistant discretizations of $\Omega$.


\subsection{Unstructured first-order system example}

\begin{figure}[t]
  \centering
  \begin{subfigure}[b]{.49\textwidth}
    \centering
  \tikzexternalenable%
  \tikzsetnextfilename{iss_tf}%
  \input{graphics/iss_tf.tikz}%
  \tikzexternaldisable%

    \caption{Transfer functions.}
    \label{fig:iss_tf}
  \end{subfigure}%
  \hfill
  \begin{subfigure}[b]{.49\textwidth}
    \centering
  \tikzexternalenable%
  \tikzsetnextfilename{iss_err}%
  \input{graphics/iss_err.tikz}%
  \tikzexternaldisable%

    \caption{Relative errors.}
    \label{fig:iss_err}
  \end{subfigure}

  \vspace{.5\baselineskip}
  \tikzexternalenable%
  \tikzsetnextfilename{iss_legend}%
  \input{graphics/iss_legend.tikz}%
  \tikzexternaldisable%

  \caption{First-order system example: All applied methods provide a similar
    approximation behavior since no special structure needs to be preserved in
    the example.
    Insignificant differences are revealed by the pointwise relative errors.}
  \label{fig:iss}
\end{figure}
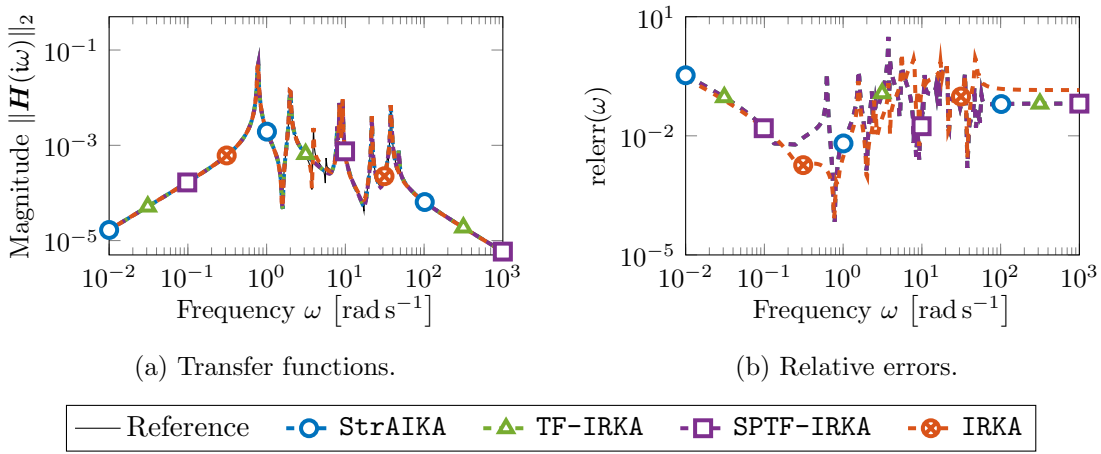

In the first example, we consider a system modeling the structural response of
the Russian Service Module of the International Space Station
(ISS)~\cite{GugAB01}.
The system is modeled using first-order differential equations~\cref{eqn:fosys}
and has $n = 270$ states, $m = 3$ inputs and $p=3$ outputs.
The transfer function is given by
\begin{equation*}
  \tf\bs = \bC \left( s \bE - \bA \right)\inv \bB.
\end{equation*}
The model is evaluated for frequencies in the range
$\convinterval = \left[\num{1e-2}, \num{1e3} \right]$\,\SI{}{\radian\per\s}.
Because of the first-order structure of the full-order model, \irka{} can be
applied in this case.
Additionally, \tfirka{}, \sptfirka{} and \straika{} are employed to compute
real-valued reduced-order models of size $r=20$ each.
Sigma plots of the transfer functions and relative approximation errors for all
models are given in \Cref{fig:iss}.
\Cref{tab:iss} summarizes the performance of the algorithms.

\begin{table}[t]
    \centering
    \caption{First-order system example:
      Comparisons of the relative, local  $\Linf$-errors,
      the required number of iterations $n_{\rm{iter}}$,
      the number of solutions of full-order $n$ linear system $n_\rm{ls}$
      and the computation time $t_{\rm{c}}$.}
    \label{tab:iss}
    \vspace{.5\baselineskip}
    
    \begin{filecontents*}{graphics/data/iss_table.csv}
algorithm,linf_error,n_iter,n_ls,t_c,mark_maxiter
\straika,0.0366925427300481,21,222,1.917002, 
\tfirka,0.0366925425283997,10,200,0.181644, 
\sptfirka,0.0366925425160824,12,132,0.265119, 
\irka,0.04217193645993,13,130,0.169957, 
\end{filecontents*}
\pgfplotstabletypeset[
        string type,
        col sep = comma,
        columns={algorithm,linf_error,n_iter,n_ls,t_c,mark_maxiter},
        every head row/.style={before row=\toprule,after row=\midrule},
        every last row/.style={after row=\bottomrule},
    ]{graphics/data/iss_table.csv}
\end{table}

All employed methods succeed in computing surrogates, which approximate the
original system up to the same degree of accuracy.
Only \irka{} obtains a different local optimum than the other methods, which
yields an insignificantly larger relative approximation error; cf.
\Cref{tab:iss}.
This meets expectations, as the first-order structure of the original system
can be represented well by the realization \tfirka{} yields.
\tfirka{} converges after only ten iterations, while the other methods require
more.
\irka{} performs the fewest decompositions of the full-order matrices and has
the shortest runtime; however, the significance of the runtime is limited for
this small example.
\straika{} requires the most iterations and therefore the most
matrix decompositions.
The runtime is longer than for the other methods.
This is a result of the additional sampling step performed in each iteration of
\straika{}, which has a measurable influence on the runtime, as $n$ is
relatively small in comparison to $r$ in this example.


\subsection{Time-delayed heated rod}

\begin{figure}[t]
	\centering
	\begin{subfigure}[b]{.49\textwidth}
		\centering
  \tikzexternalenable%
  \tikzsetnextfilename{delay_tf}%
  \input{graphics/delay_tf.tikz}%
  \tikzexternaldisable%

		\caption{Transfer functions.}
		\label{fig:delay_tf}
	\end{subfigure}%
	\hfill
	\begin{subfigure}[b]{.49\textwidth}
		\centering
  \tikzexternalenable%
  \tikzsetnextfilename{delay_err}%
  \input{graphics/delay_err.tikz}%
  \tikzexternaldisable%

		\caption{Relative errors.}
		\label{fig:delay_err}
	\end{subfigure}
	
	\vspace{.5\baselineskip}
  \tikzexternalenable%
  \tikzsetnextfilename{delay_legend}%
  \input{graphics/delay_legend.tikz}%
  \tikzexternaldisable%

	\caption{Time-delay system example:
      The unstructured approximation computed by \tfirka{} is around two orders
      of magnitude less accurate than the structure-preserving reduced-order
      models computed by \straika{} and \sptfirka{} for smaller to medium
      frequencies.
      Overall \straika{} provides the most accurate approximation.}
	\label{fig:delay}
\end{figure}
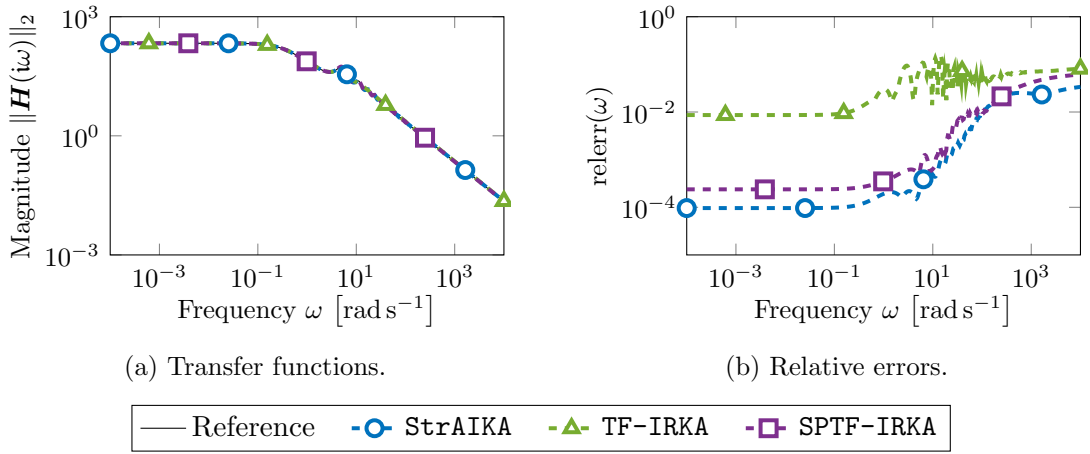

Here, we consider a model of a heated rod with distributed control and
homogeneous Dirichlet boundary conditions, which is cooled by delayed feedback.
This system has also been analyzed in~\cite{BenGW22}.
A discretization of the underlying partial differential equation leads to the
transfer function
\begin{equation*}
  \tf\bs = \bC \left( s \bE - \bA_{0} - \operatorname{e}^{-\tau s}
    \bA_{\rm{d}} \right)\inv \bB,
\end{equation*}
with $n = 1\,000\,000$ states, $m = 5$ inputs and $p = 4$ outputs.
The time delay is $\tau = 1$ in this example.
The input matrix $\bB$ has a block structure such that the rod is heated
uniformly at different sections by the inputs.
The outputs are the average temperatures on these sections.
For this example, the frequency range
$\convinterval = \left[\num{1e-4}, \num{1e4}\right]$\,\SI{}{\radian\per\s} is
considered.
For the experiments, we fix the reduced order to $r = 10$ and compute
reduced-order models with \straika{}, \tfirka{} and \sptfirka{}.
The system cannot be transformed into an equivalent system with first-order
structure, therefore \irka{} cannot be applied in this case.
The initial interpolation points are distributed linearly equidistant in
$\qim \left[\num{1e-4}, \num{1e2}\right]$.
\Cref{fig:delay} plots the maximum singular values of the transfer functions
and the error systems.
Further results are given in \Cref{tab:delay}.

\begin{table}[t]
  \centering
  \caption{Time-delay system example:
    Comparisons of the relative, local  $\Linf$-errors,
    the required number of iterations $n_{\rm{iter}}$,
    the number of solutions of full-order $n$ linear system $n_\rm{ls}$
    and the computation time $t_{\rm{c}}$.}
  \label{tab:delay}
  \vspace{.5\baselineskip}
    
    \begin{filecontents*}{graphics/data/delay_table.csv}
algorithm,linf_error,n_iter,n_ls,t_c,mark_maxiter
\straika,0.000101105936687689,37,190,247.101835, 
\tfirka,0.020574550747899,73,730,367.244314, 
\sptfirka,0.000238103757728217,3,29,37.562543, 
\end{filecontents*}
\pgfplotstabletypeset[
        string type,
        col sep = comma,
        columns={algorithm,linf_error,n_iter,n_ls,t_c,mark_maxiter},
        every head row/.style={before row=\toprule,after row=\midrule},
        every last row/.style={after row=\bottomrule},
    ]{graphics/data/delay_table.csv}
\end{table}

In this example, \straika{} and \sptfirka{} compute reduced-order models with
comparable accuracy.
However, \straika{} provides the smallest worst case error by around a factor
of two as shown in \Cref{tab:delay}.
As expected, the first-order realization computed by \tfirka{} cannot capture
the dynamics of the delay system well and, therefore, provides an approximation
that is around two orders of magnitude less accurate than the models computed
by \straika{} and \sptfirka{}.
Because of its rapid convergence, the runtime of \sptfirka{} is considerably
lower than for the other two algorithms.
But also the runtime of \straika{} is significantly smaller than for \tfirka{}.
Note, that most of the system dynamics happen in the considered frequency
range.
Therefore, \straika{} quickly reaches the maximum reduced order and has to
choose the $10$ most dominant poles out of approximately $150$ eigenvalues of
the Loewner realization in each iteration.
This additional effort directly affects the runtime of \straika{}.
Additionally, the convergence is slow compared to \sptfirka{}.


\subsection{Viscoelastic beam}

\begin{figure}[t]
	\centering
	\begin{subfigure}[b]{.49\textwidth}
		\centering
  \tikzexternalenable%
  \tikzsetnextfilename{visco_tf}%
  \input{graphics/visco_tf.tikz}%
  \tikzexternaldisable%

		\caption{Transfer functions.}
		\label{fig:visco_tf}
	\end{subfigure}%
	\hfill
	\begin{subfigure}[b]{.49\textwidth}
		\centering
  \tikzexternalenable%
  \tikzsetnextfilename{visco_err}%
  \input{graphics/visco_err.tikz}%
  \tikzexternaldisable%

		\caption{Relative errors.}
		\label{fig:visco_err}
	\end{subfigure}
	
	\vspace{.5\baselineskip}
  \tikzexternalenable%
  \tikzsetnextfilename{visco_legend}%
  \input{graphics/visco_legend.tikz}%
  \tikzexternaldisable%

  \caption{Viscoelastic beam example:
    All methods succeed in computing reasonably accurate reduced-order models.
    The relative approximation error of the reduced-order model obtained by
    \straika{} is around three orders of magnitude smaller compared to the other
    methods.}
  \label{fig:visco}
\end{figure}
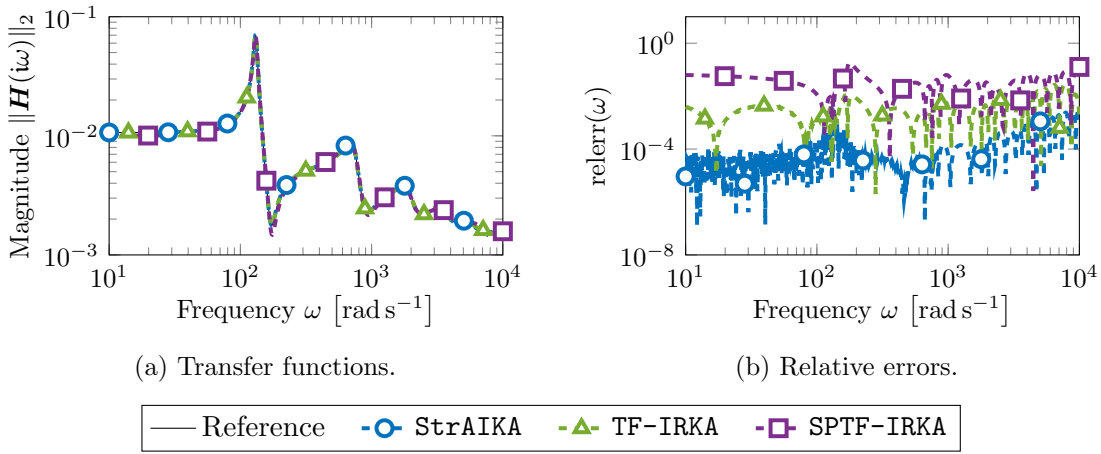

This example models a flexible beam with viscoelastic core.
The beam of length $l = \SI{0.21}{\meter}$ has a symmetric sandwich structure
consisting of two layers of cold rolled steel surrounding a viscoelastic
ethylene-propylene-diene core~\cite{VanMM13}; the beam is clamped at one side.
After discretization, the transfer function of the system is given by
\begin{equation*}
  \tf\bs = \bC \left( s^{2} \bM + \bK + \frac{G_{0} + G_{\infty}
    \left(s\tau\right)^{\alpha}}
    {1 + \left(s \tau \right)^{\alpha}} \bG \right)\inv \bB.
\end{equation*}
The model has $n = 3\,360$ states, $m = 1$ input and $p = 1$ output.
The beam is excited by a single load at its free end and the displacement is
measured at the same location, resulting in output and input mappings
$\bC = 100 \cdot \bB\trans = \begin{bmatrix} 0 & \cdots & 0 & 1 \end{bmatrix}$.
The matrices $\bM$, $\bK$, $\bG$ are available from~\cite{HigNT19}.
In this example, we limit the frequency range of interest to
$\convinterval = \left[\num{10}, \num{1e4}\right]$\,\SI{}{\radian\per\s}.
Note, that the system has poles, which lie outside of this range.
No maximum reduced order $r$ is set in this case, so \straika{} determines it
in an adaptive way.
The initial interpolation points are a single complex conjugate pair, where the
absolute value of its imaginary part is located in the middle of
$\convinterval$.
The automatically determined order is used for the experiments with \tfirka{}
and \sptfirka{}, where the $\lceil r/2 \rceil$ initial expansion points and
their complex conjugates are distributed logarithmically equidistant in
$\convinterval$.
The sigma plots for the frequency responses of the reference and the
reduced-order models as well as the corresponding errors are given in
\Cref{fig:visco}.
The performance of the methods is shown in \Cref{tab:visco}.

\begin{table}[t]
  \centering
  \caption{Viscoelastic beam example:
    Comparisons of the relative, local  $\Linf$-errors,
    the required number of iterations $n_{\rm{iter}}$,
    the number of solutions of full-order $n$ linear system $n_\rm{ls}$
    and the computation time $t_{\rm{c}}$.
    The $\ast$ marks experiments, where the maximum number of
    iterations has been reached without convergence.}
  \label{tab:visco}
  \vspace{.5\baselineskip}

    \begin{filecontents*}{graphics/data/sandwich_table.csv}
algorithm,linf_error,n_iter,n_ls,t_c,mark_maxiter
\straika,0.000913874693057171,10,74,0.709724, 
\tfirka,0.00103672128471925,18,288,0.970171, 
\sptfirka,0.00959688757116932,50,457,8.276918,$\ast$
\end{filecontents*}
\pgfplotstabletypeset[
    	string type,
    	col sep = comma,
    	columns={algorithm,linf_error,n_iter,n_ls,t_c,mark_maxiter},
    	every head row/.style={before row=\toprule,after row=\midrule},
    	every last row/.style={after row=\bottomrule},
      ]{graphics/data/sandwich_table.csv}
\end{table}

\straika{} converges after ten iterations to a model of size $r = 15$, i.e.,
the reduced-order models computed by \tfirka{} and \sptfirka{} have order
$r = 16$.
All three algorithms produce reasonably accurate models regarding the
reference, while the model computed by \straika{} is around three orders of
magnitude more accurate for most frequencies; cf. \Cref{fig:visco_err}.
Also the worst case approximation error of \straika{} is around one order
of magnitude better than \sptfirka{}; see \Cref{tab:visco}.
This can be explained by the fact that \straika{} places all interpolation
points in $\convinterval$, while some of the interpolation points obtained by
the other two algorithms are located in the frequency region above
\SI{1e4}{\radian\per\s}, leading to a higher error in $\convinterval$.
\sptfirka{} did not converge after $50$ iterations, resulting in a high
computation time.
In terms of computational effort, \straika{} is clearly advantageous in this
example.


\subsection{Radio frequency gun}

\begin{figure}[t]
\centering
\begin{subfigure}[b]{.49\textwidth}
	\centering
  \tikzexternalenable%
  \tikzsetnextfilename{gun_tf}%
  \input{graphics/gun_tf.tikz}%
  \tikzexternaldisable%

	\caption{Transfer functions.}
	\label{fig:gun_tf}
\end{subfigure}%
\hfill
\begin{subfigure}[b]{.49\textwidth}
	\centering
  \tikzexternalenable%
  \tikzsetnextfilename{gun_err}%
  \input{graphics/gun_err.tikz}%
  \tikzexternaldisable%

	\caption{Relative errors.}
	\label{fig:gun_err}
\end{subfigure}

\vspace{.5\baselineskip}
  \tikzexternalenable%
  \tikzsetnextfilename{gun_legend}%
  \input{graphics/gun_legend.tikz}%
  \tikzexternaldisable%

  \caption{Radio frequency gun example:
    The reduced-order model computed by \straika{} shows a relatively uniform
    accuracy in $\convinterval$, while the surrogate obtained form \sptfirka{}
    is only accurate for less than \SI{60}{\radian\per\s}.
    The first-order system computed by \tfirka{} is not able to approximate the
    dynamics of the original system at all.}
  \label{fig:gun}
\end{figure}
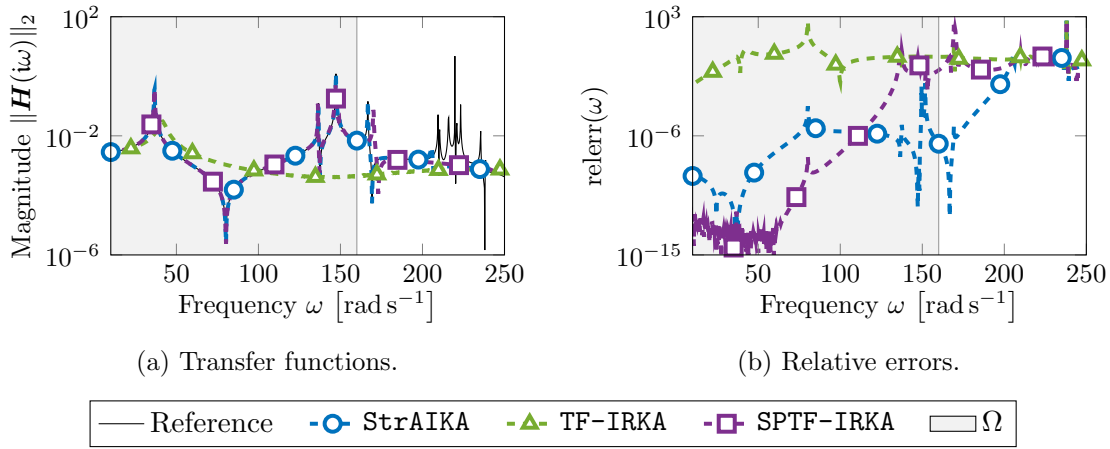

As the last example, we consider a radio frequency gun as described
in~\cite{Lia07}.
Discretizing the system leads to a transfer function
\begin{equation*}
  \tf\bs = \bC \left( s^{2} \bM + \qim \left(s^{2} -
    \sigma_{1}^{2} \right)^{\frac{1}{2}} \bW_{1} + \qim\left(s^{2} -
    \sigma_{2}^{2} \right)^{\frac{1}{2}} \bW_{2} + \bK \right)\inv \bB,
\end{equation*}
where $\sigma_{1} = 0.0$ and $\sigma_{2} = 108.8774$.
The full-order model has $n = 9\,956$ states and the matrices
$\bM, \bW_{1}, \bW_{2}, \bK$ are taken from~\cite{HigNT19}.
We consider $m = 1$ input and $p = 1$ output.
The input vector is populated with values drawn from the standard normal
distribution and the system output is measured at the first degree of freedom,
corresponding to $\bC = \begin{bmatrix} 1 & 0 & \ldots & 0 \end{bmatrix}$.
These vectors are also part of the code package supplementing this
article~\cite{supAumW23}.
The reduced-order models are computed to approximate the original system in
$\convinterval = \left[1, 160\right]$\,\SI{}{\radian\per\s} and the necessary
order $r$ for the surrogate is automatically determined by \straika{}.
One initial complex conjugate pair of interpolation points is placed in the
middle of $\convinterval$.
The initial $\lceil r/2 \rceil$ pairs of interpolation points for \tfirka{} and
\sptfirka{} are distributed linearly equidistant in $\convinterval$.
The sigma plots for the frequency responses of the reference and the
reduced-order models as well as the corresponding relative approximation errors
are shown in \Cref{fig:gun}.
The local, relative $\Linf$-errors and computational costs are presented in
\Cref{tab:gun}.

\begin{table}[t]
    \centering
    \caption{Radio frequency gun example:
    Comparisons of the relative, local  $\Linf$-errors,
    the required number of iterations $n_{\rm{iter}}$,
    the number of solutions of full-order $n$ linear system $n_\rm{ls}$
    and the computation time $t_{\rm{c}}$.
    The $\ast$ marks experiments, where the maximum number of
    iterations has been reached without convergence.}
    \label{tab:gun}
    \vspace{.5\baselineskip}

    \begin{filecontents*}{graphics/data/gun_table.csv}
algorithm,linf_error,n_iter,n_ls,t_c,mark_maxiter
\straika,0.000173015039226412,50,416,123.327627,$\ast$
\tfirka,0.999844667946225,100,3000,2909.537882,$\ast$
\sptfirka,0.999775790485055,50,1316,392.115224,$\ast$
\end{filecontents*}
\pgfplotstabletypeset[
    	string type,
    	col sep = comma,
    	columns={algorithm,linf_error,n_iter,n_ls,t_c,mark_maxiter},
    	every head row/.style={before row=\toprule,after row=\midrule},
    	every last row/.style={after row=\bottomrule},
      ]{graphics/data/gun_table.csv}
\end{table}

\straika{} automatically determines a reduced order of $r = 30$ and computes a
reduced-order model, which is more accurate in the frequency range of interest
than the reduced-order models computed by the other two algorithms.
Placing the interpolation points only inside $\convinterval$ leads to a higher
accuracy in this region, while the approximation deteriorates for higher
frequencies.
Although \straika{} does not converge after 50 iterations, the reduced-order
model is more than reasonably accurate, with its worst case approximation error
being four orders of magnitude smaller than for the other methods; see
\Cref{tab:gun}.
The other two methods also do not converge in the allowed number of iterations.
This example clearly shows the increased computational effort required for
\tfirka{} as well as \sptfirka{}.
In particular, \sptfirka{} takes more than three times longer than \straika{}
for the same number of iterations.
This is a result of the large number of iterations required for the inner
\tfirka{} inside of \sptfirka{} leading to a large number of linear solves to
be performed.
The unstructured \tfirka{} fails in computing a surrogate that approximates the
transfer function of the full-order model.
The reduced-order model computed by \sptfirka{} approximates the
full-order model well in the lower frequency region, but the overall accuracy in
$\convinterval$ is lower than observed for the model computed by \straika{}.
The comparison to \tfirka{} again shows the importance of
structure-preserving model order reduction strategies.


\section{Conclusions}
\label{sec:conclusions}

In this paper, we proposed \straika{}, a new algorithm that computes
structure-preserving reduced-order models of systems with arbitrary transfer
function structure in an iterative way.
Similar to \irka{}-like methods, \straika{} uses mirror images of intermediate
reduced-order models as interpolation points for the next iteration.
In each iteration, the Loewner framework is used to compute first-order
realizations of transfer function data collected from the current reduced-order
model, whose eigenvalues are the basis for the interpolation points in the next
iteration.
\straika{} automatically determines a size for the reduced-order model by
considering the eigenvalues located in a given frequency range of interest.
This ensures a reasonable approximation of the system dynamics at least in
these regions.
\straika{} is completely agnostic to the actual transfer function structure of
the original model and does not require any transfer function derivatives.

We demonstrated the versatility and effectiveness of \straika{} in four
numerical examples with different internal structures.
The benchmark systems model structural vibration, heat transfer with internal
delay, viscoelasticity, and radio wave propagation.
\straika{} showed comparable or even significantly better accuracy with respect
to established \irka{}-like methods.
Especially, if only a limited frequency range is of interest for the
application, \straika{} easily outperforms methods that optimize the
approximation error under the $\Htwo$-norm both in terms of accuracy and
required computational effort.

An open question that has not been covered in this paper is the preservation of
additional system properties like stability.
However, the preservation of such system properties has only been solved for
particular system structures such as first-order systems or second-order
systems with mechanical matrix structure.
For the case of generally structured systems as considered in this paper, no
solution to this problem for projection-based model reduction is known yet.
Another idea for future investigations is the reduction of the computational
costs of \straika{} by considering an additional layer of approximation as
it has been done for \irka{}-like methods in~\cite{AumM22, CasL18}.


\addcontentsline{toc}{section}{References}
\bibliographystyle{plainurl}
\bibliography{bibtex/myref}

\end{document}

%% file: graphics/region_spectrum.tikz
\begin{tikzpicture}[font = \small]
  \pgfplotstableread[col sep=comma]{graphics/data/region_ev.csv}{\dataev}
  \pgfplotstableread[col sep=comma]{graphics/data/region_s0.csv}{\datas}

  \begin{axis}[name=evs,
    height = .2\textheight,
    width = .94\textwidth,
    xlabel = {\qreal{\lambda}},
    ylabel = {\qimag{\lambda}},
    ylabel style = {yshift = -.4em},
    xlabel style = {yshift = .4em},
	xmin = -100,
    xmax = 100,
    ymin = -5e3,
    ymax = 5e3,
    axis on top,
  ]

    \addplot[only marks, mark = o, black, mark options = {scale = 1.5}]
      table[x = real, y = imag] {\dataev};

    \addplot[only marks, mark = o, color = colD, mark options = {scale=1.5}]
      table[x = real, y = imag] {\datas};
    
    \addplot[only marks, mark = asterisk, color= colD, mark options={scale=1.5}]
      table[x expr = -1*\thisrow{real}, y = imag] {\datas};

    \addplot[name path = min, mark = none, black!40, forget plot]
      coordinates {
        (\pgfkeysvalueof{/pgfplots/xmin}, 4000)
        (\pgfkeysvalueof{/pgfplots/xmax}, 4000)
      };%
    \addplot [name path = max, mark = none, black!40, forget plot]
      coordinates {
        (\pgfkeysvalueof{/pgfplots/xmin}, 3000)
        (\pgfkeysvalueof{/pgfplots/xmax}, 3000)
      };%
	\addplot[fill = black!10, area legend] fill between [of = min and max];
 
    \addplot[name path = minb, mark = none, black!40, forget plot]
      coordinates {
        (\pgfkeysvalueof{/pgfplots/xmin}, -4000)
        (\pgfkeysvalueof{/pgfplots/xmax}, -4000)
      };%
    \addplot[name path = maxb, mark = none, black!40, forget plot]
      coordinates {
        (\pgfkeysvalueof{/pgfplots/xmin}, -3000)
        (\pgfkeysvalueof{/pgfplots/xmax}, -3000)
      };%
	\addplot[fill = black!5, area legend] fill between [of = minb and maxb];
    
  \end{axis}
\end{tikzpicture}

%% file: graphics/region_tf.tikz
\begin{tikzpicture}[font = \small]
  \pgfplotstableread[col sep=comma]{graphics/data/region.csv}{\data}

  \begin{axis}[
    xmode=log, ymode=log,
    height = .2\textheight,
    width  = .92\textwidth,
    xlabel = {Frequency},
    x unit = \radian\per\s,
    ylabel = {Magnitude},
    ylabel style = {yshift = -.4em},
    xlabel style = {yshift = .4em},
	xmin = 1e2,
    xmax = 1e4,
    ymin = 1e-10,
    ymax = 1e-5,
	axis on top,
	cycle list name  = qcolorlist,
	cycle list shift = -1,
    legend cell align = left,
    max space between ticks = 30
  ]
	
    \addplot+[no marks, color=black] table[x=s, y=res] {\data};
    \addplot+[dashed, no marks, color=colA] table[x=s, y=res_r] {\data};

    \addplot[only marks, mark=o, black, mark options={scale=1.5}] coordinates {(1,1)};

    \addplot[only marks, mark=asterisk, color=colD, mark options={scale=1.5}] coordinates {(1,1)};

    \addplot[only marks, mark=x, draw=black!80, fill=colB!40, mark options={scale=1.5}] coordinates {(1,1)};
    
    \plotinterval{3000}{4000}
\end{axis}
\end{tikzpicture}

%% file: graphics/region_err.tikz
\begin{tikzpicture}[font = \small]
  \pgfplotstableread[col sep=comma]{graphics/data/region.csv}{\data}
  \pgfplotstableread[col sep=comma]{graphics/data/region_vals.csv}{\datasnew}

  \begin{axis}[
    xmode=log, ymode=log,
    height = .2\textheight,
    width  = .92\textwidth,
    xlabel = {Frequency},
    x unit = \radian\per\s,
    ylabel = {Magntiude},
    ylabel style = {yshift = -.4em},
    xlabel style = {yshift = .4em},
    xmin = 1e2,
    xmax = 1e4,
    ymin = 1e-12,
    ymax = 1e4,
    axis on top,
    legend pos = outer north east, legend cell align = left,
    cycle list name=qcolorlist,
  ]
  
  \addplot+[dashed, no marks] table[x=s, y=err] {\data};
  
  \addplot[only marks, mark=x, draw=black!80, fill=colB!40, mark options={scale=1.5}] table[x=s, y=err] {\datasnew};
  
  \plotinterval{3000}{4000}
\end{axis}
\end{tikzpicture}

%% file: graphics/region_legend.tikz
\begin{tikzpicture}
  \begin{axis}[%
    hide axis,
    width  = 1mm,
    height = 1mm,
    scale only axis,
    xmin = 0,
    xmax = 1,
    ymin = 0,
    ymax = 1,
    legend columns = -1, 
    legend style   = {
      at     = {(0,0)},
      anchor = center,
      /tikz/every even column/.append style = {column sep = 1em}},
    legend cell align  = {left},
    cycle list name    = qcolorlist,
	cycle list shift   = -1,
    clip mode          = individual]
    
    \addplot+[no marks, color=black] coordinates {(0,0)};
    \addlegendentry{Reference}

    \addplot+[dashed, no marks, color=colA] coordinates {(0,0)};
    \addlegendentry{ROM}

    \addplot[only marks, mark = o, black, mark options = {scale = 1.5}]
      coordinates {(0,0)};
    \addlegendentry{$\lambda_{j}$}

    \addplot[only marks, mark = asterisk, color=colD, mark options={scale=1.5}]
      coordinates {(0,0)};
    \addlegendentry{$\sigma_{j}$}

    \addplot[only marks, mark = x, draw = black!80, fill=colB!40,
      mark options={scale=1.5}] coordinates {(0,0)};
    \addlegendentry{$\qimag{\sigma_{j}}$}

    \plotinterval{0}{0}
    \addlegendentry{\convinterval}
  \end{axis}
\end{tikzpicture}

%% file: graphics/iss_tf.tikz
\begin{tikzpicture}[font = \small]
  \pgfplotstableread[col sep=comma]{graphics/data/iss_sigma.csv}{\data}

  \begin{axis}[
    height = .2\textheight,
    width  = .92\textwidth,
    xmode = log,
    ymode = log,
	xlabel = {Frequency $\omega$},
	x unit = {\radian\per\s},
	ylabel = {Magnitude $\sigmaplot$},
    ylabel style = {yshift = -.4em},
    xlabel style = {yshift = .4em},
	xmin=1e-2, xmax=1e3,
	ymin=5e-6, ymax=5e-1,
		cycle list name  =qcolorlist,
		cycle list shift=-1,
		max space between ticks=30,
    mark repeat={200}
  ]
		\addplot[mark=none, color=black] table[x=s,  y=res] {\data};
		\addplot+[dashed] table[x=s, y=res_r] {\data};
		\addplot+[dashed, mark phase=50, dash phase=1pt] table[x=s, y=res_tfi] {\data};
		\addplot+[dashed, mark phase=100, dash phase=2pt] table[x=s, y=res_sptfi] {\data};
		\addplot+[dashed, mark phase=150, dash phase=3pt] table[x=s, y=res_i] {\data};
	\end{axis}
\end{tikzpicture}
 

%% file: graphics/iss_err.tikz
\begin{tikzpicture}[font = \small]
  \pgfplotstableread[col sep=comma]{graphics/data/iss_sigma.csv}{\data}
  
  \begin{axis}[
    height = .2\textheight,
    width  = .92\textwidth,
    xmode=log,
    ymode=log,
	xlabel = {Frequency $\omega$},
	x unit = {\radian\per\s},
	ylabel = {$\sigmaerr$},
    ylabel style = {yshift = -.4em},
    xlabel style = {yshift = .4em},
	xmin=1e-2, xmax=1e3,
	ymin=1e-5, ymax=1e1,
	legend cell align = left,
	legend columns=-1,
	cycle list name=qcolorlist,
	max space between ticks=30,
	mark repeat={200}
  ]
    \addplot+[dashed] table[x=s, y=err] {\data};
	\addplot+[dashed, mark phase=50] table[x=s, y=err_tfi] {\data};
	\addplot+[dashed, mark phase=100] table[x=s, y=err_sptfi] {\data};
	\addplot+[dashed, mark phase=150] table[x=s, y=err_i] {\data};
  \end{axis}
\end{tikzpicture}

%% file: graphics/iss_legend.tikz
\begin{tikzpicture}
  \begin{axis}[%
    hide axis,
    width  = 1mm,
    height = 1mm,
    scale only axis,
    xmin = 0,
    xmax = 1,
    ymin = 0,
    ymax = 1,
    legend columns = -1, 
    legend style   = {
      at     = {(0,0)},
      anchor = center,
      /tikz/every even column/.append style = {column sep = 1em}},
    legend cell align  = {left},
    cycle list name    = qcolorlist,
	cycle list shift   = -1,
    clip mode          = individual]

    \addplot[mark=none, color=black] coordinates {(0,0)};
    \addlegendentry{Reference}
    
    \pgfplotsinvokeforeach{1, ..., 4}{\addplot+[dashed] coordinates {(0,0)};}
    \addlegendentry{\straika}
    \addlegendentry{\tfirka}
    \addlegendentry{\sptfirka}
    \addlegendentry{\irka}
  \end{axis}
\end{tikzpicture}

%% file: graphics/delay_tf.tikz
\begin{tikzpicture}[font = \small]
  \pgfplotstableread[col sep=comma]{graphics/data/delay_sigma.csv}{\data}

  \begin{axis}[
    height = .2\textheight,
    width  = .92\textwidth,
    xmode = log,
    ymode = log,
	xlabel = {Frequency $\omega$},
	x unit = {\radian\per\s},
	ylabel = {Magnitude $\sigmaplot$},
    ylabel style = {yshift = -.4em},
    xlabel style = {yshift = .4em},
    xmin=1e-4, xmax=1e4,
    ymin=1e-3, ymax=1e3,
		cycle list name  =qcolorlist,
		cycle list shift=-1,
		max space between ticks=30,
    mark repeat={150}
  ]
		\addplot[mark=none, color=black] table[x=s,  y=res] {\data};
		\addplot+[dashed] table[x=s, y=res_r] {\data};
		\addplot+[dashed, mark phase=50, dash phase=1pt] table[x=s, y=res_tfi] {\data};
		\addplot+[dashed, mark phase=100, dash phase=2pt] table[x=s, y=res_sptfi] {\data};
	\end{axis}
\end{tikzpicture}
 

%% file: graphics/delay_err.tikz
\begin{tikzpicture}[font = \small]
  \pgfplotstableread[col sep=comma]{graphics/data/delay_sigma.csv}{\data}
  
  \begin{axis}[
    height = .2\textheight,
    width  = .92\textwidth,
    xmode=log,
    ymode=log,
	xlabel = {Frequency $\omega$},
	x unit = {\radian\per\s},
	ylabel = {$\sigmaerr$},
    ylabel style = {yshift = -.4em},
    xlabel style = {yshift = .4em},
    xmin=1e-4, xmax=1e4,
    ymin=1e-5, ymax=1,
	legend cell align = left,
	legend columns=-1,
	cycle list name=qcolorlist,
	max space between ticks=30,
	mark repeat={150}
  ]
    \addplot+[dashed] table[x=s, y=err] {\data};
	\addplot+[dashed, mark phase=50] table[x=s, y=err_tfi] {\data};
	\addplot+[dashed, mark phase=100] table[x=s, y=err_sptfi] {\data};
  \end{axis}
\end{tikzpicture}

%% file: graphics/delay_legend.tikz
\begin{tikzpicture}
  \begin{axis}[%
    hide axis,
    width  = 1mm,
    height = 1mm,
    scale only axis,
    xmin = 0,
    xmax = 1,
    ymin = 0,
    ymax = 1,
    legend columns = -1, 
    legend style   = {
      at     = {(0,0)},
      anchor = center,
      /tikz/every even column/.append style = {column sep = 1em}},
    legend cell align  = {left},
    cycle list name    = qcolorlist,
	cycle list shift   = -1,
    clip mode          = individual]

    \addplot[mark=none, color=black] coordinates {(0,0)};
    \addlegendentry{Reference}
    
    \pgfplotsinvokeforeach{1, ..., 3}{\addplot+[dashed] coordinates {(0,0)};}
    \addlegendentry{\straika}
    \addlegendentry{\tfirka}
    \addlegendentry{\sptfirka}
  \end{axis}
\end{tikzpicture}

%% file: graphics/visco_tf.tikz
\begin{tikzpicture}[font = \small]
  \pgfplotstableread[col sep=comma]{graphics/data/sandwich_sigma.csv}{\data}

  \begin{axis}[
    height = .2\textheight,
    width  = .92\textwidth,
    xmode = log,
    ymode = log,
	xlabel = {Frequency $\omega$},
	x unit = {\radian\per\s},
	ylabel = {Magnitude $\sigmaplot$},
    ylabel style = {yshift = -.4em},
    xlabel style = {yshift = .4em},
    xmin=1e1, xmax=1e4,
    ymin=1e-3, ymax=1e-1,
		cycle list name  =qcolorlist,
		cycle list shift=-1,
		max space between ticks=30,
    mark repeat={150}
  ]
		\addplot[mark=none, color=black] table[x=s,  y=res] {\data};
		\addplot+[dashed] table[x=s, y=res_r] {\data};
		\addplot+[dashed, mark phase=50, dash phase=1pt] table[x=s, y=res_tfi] {\data};
		\addplot+[dashed, mark phase=100, dash phase=2pt] table[x=s, y=res_sptfi] {\data};
	\end{axis}
\end{tikzpicture}
 

%% file: graphics/visco_err.tikz
\begin{tikzpicture}[font = \small]
  \pgfplotstableread[col sep=comma]{graphics/data/sandwich_sigma.csv}{\data}
  
  \begin{axis}[
    height = .2\textheight,
    width  = .92\textwidth,
    xmode=log,
    ymode=log,
	xlabel = {Frequency $\omega$},
	x unit = {\radian\per\s},
	ylabel = {$\sigmaerr$},
    ylabel style = {yshift = -.4em},
    xlabel style = {yshift = .4em},
    xmin=1e1, xmax=1e4,
    ymin=1e-8, ymax=1e1,
	legend cell align = left,
	legend columns=-1,
	cycle list name=qcolorlist,
	max space between ticks=30,
	mark repeat={150}
  ]
    \addplot+[dashed] table[x=s, y=err] {\data};
	\addplot+[dashed, mark phase=50] table[x=s, y=err_tfi] {\data};
	\addplot+[dashed, mark phase=100] table[x=s, y=err_sptfi] {\data};
  \end{axis}
\end{tikzpicture}

%% file: graphics/visco_legend.tikz
\begin{tikzpicture}
  \begin{axis}[%
    hide axis,
    width  = 1mm,
    height = 1mm,
    scale only axis,
    xmin = 0,
    xmax = 1,
    ymin = 0,
    ymax = 1,
    legend columns = -1, 
    legend style   = {
      at     = {(0,0)},
      anchor = center,
      /tikz/every even column/.append style = {column sep = 1em}},
    legend cell align  = {left},
    cycle list name    = qcolorlist,
	cycle list shift   = -1,
    clip mode          = individual]

    \addplot[mark=none, color=black] coordinates {(0,0)};
    \addlegendentry{Reference}
    
    \pgfplotsinvokeforeach{1, ..., 3}{\addplot+[dashed] coordinates {(0,0)};}
    \addlegendentry{\straika}
    \addlegendentry{\tfirka}
    \addlegendentry{\sptfirka}
  \end{axis}
\end{tikzpicture}

%% file: graphics/gun_tf.tikz
\begin{tikzpicture}[font = \small]
  \pgfplotstableread[col sep=comma]{graphics/data/gun_sigma.csv}{\data}

  \begin{axis}[
    height = .2\textheight,
    width  = .92\textwidth,
    ymode = log,
	xlabel = {Frequency $\omega$},
	x unit = {\radian\per\s},
	ylabel = {Magnitude $\sigmaplot$},
    ylabel style = {yshift = -.4em},
    xlabel style = {yshift = .4em},
    xmin=1e1, xmax=250,
    ymin=1e-6, ymax=1e2,
	axis on top,
    cycle list name  =qcolorlist,
    cycle list shift=-1,
    max space between ticks=30,
    mark repeat={150}
  ]
		\addplot[mark=none, color=black] table[x=s,  y=res] {\data};
		\addplot+[dashed] table[x=s, y=res_r] {\data};
		\addplot+[dashed, mark phase=50, dash phase=1pt] table[x=s, y=res_tfi] {\data};
		\addplot+[dashed, mark phase=100, dash phase=2pt] table[x=s, y=res_sptfi] {\data};
        \plotinterval{1}{160}
	\end{axis}
\end{tikzpicture}
 

%% file: graphics/gun_err.tikz
\begin{tikzpicture}[font = \small]
  \pgfplotstableread[col sep=comma]{graphics/data/gun_sigma.csv}{\data}
  
  \begin{axis}[
    height = .2\textheight,
    width  = .92\textwidth,
    ymode=log,
	xlabel = {Frequency $\omega$},
	x unit = {\radian\per\s},
	ylabel = {$\sigmaerr$},
    ylabel style = {yshift = -.4em},
    xlabel style = {yshift = .4em},
    xmin=1e1, xmax=250,
    ymin=1e-15, ymax=1e3,
    axis on top,
	cycle list name=qcolorlist,
	max space between ticks=30,
	mark repeat={150}
  ]
    \addplot+[dashed] table[x=s, y=err] {\data};
	\addplot+[dashed, mark phase=50] table[x=s, y=err_tfi] {\data};
	\addplot+[dashed, mark phase=100] table[x=s, y=err_sptfi] {\data};
    \plotinterval{1}{160}
  \end{axis}
\end{tikzpicture}

%% file: graphics/gun_legend.tikz
\begin{tikzpicture}
  \begin{axis}[%
    hide axis,
    width  = 1mm,
    height = 1mm,
    scale only axis,
    xmin = 0,
    xmax = 1,
    ymin = 0,
    ymax = 1,
    legend columns = -1, 
    legend style   = {
      at     = {(0,0)},
      anchor = center,
      /tikz/every even column/.append style = {column sep = 1em}},
    legend cell align  = {left},
    cycle list name    = qcolorlist,
	cycle list shift   = -1,
    clip mode          = individual]

    \addplot[mark=none, color=black] coordinates {(0,0)};
    \addlegendentry{Reference}
    
    \pgfplotsinvokeforeach{1, ..., 3}{\addplot+[dashed] coordinates {(0,0)};}
    \addlegendentry{\straika}
    \addlegendentry{\tfirka}
    \addlegendentry{\sptfirka}

    \plotinterval{0}{0}
    \addlegendentry{\convinterval}
  \end{axis}
\end{tikzpicture}